\documentclass{IEEEtran4PSCC}
\usepackage[T1]{fontenc}
\usepackage[utf8]{inputenc}
\usepackage{graphicx}
\usepackage{amsmath,amssymb,amsfonts,amsthm,bm,mathtools,dsfont}
\usepackage{bigints,latexsym}
\usepackage[nointegrals]{wasysym}
\usepackage[dvipsnames]{xcolor}
\usepackage[font=small,labelfont=bf]{caption}
\usepackage{subcaption}
\usepackage{multirow,booktabs}
\usepackage{algorithm,algpseudocode}
\usepackage{enumerate}
\usepackage{tcolorbox}
\usepackage{xurl}
\usepackage{varioref}
\usepackage[hidelinks]{hyperref}
\usepackage{nameref}
\usepackage[nameinlink,noabbrev]{cleveref}

\definecolor{mustard}{HTML}{707238}
\definecolor{PaleGreen}{HTML}{98FB98}
\definecolor{LightRed}{HTML}{FE7A7A}
\definecolor{LightYellow}{HTML}{EFFA86}
\definecolor{YellowOrange}{HTML}{FFAE42}
\definecolor{OliveGreen}{HTML}{408000}
\definecolor{dkgreen}{rgb}{0,0.6,0}
\definecolor{gray}{rgb}{0.5,0.5,0.5}
\definecolor{mauve}{rgb}{0.58,0,0.82}
\definecolor{darkred}{RGB}{139,0,0}
\let\oldReturn\Return
\renewcommand{\Return}{\State\oldReturn}

\title{Real-Time Dynamic $N-1$ Screening: Identifying High-Risk Lines and Transformers After Common Faults}
\author{
\IEEEauthorblockN{Ayrton Almada \IEEEauthorrefmark{1}\IEEEauthorrefmark{2}, Laurent Pagnier \IEEEauthorrefmark{1}, Igal Goldshtein \IEEEauthorrefmark{3}, Saif R. Kazi \IEEEauthorrefmark{2},
and Michael (Misha) Chertkov \IEEEauthorrefmark{1}}
\IEEEauthorblockA{\IEEEauthorrefmark{1} Program in Applied Mathematics and Department of Mathematics,\\
University of Arizona, Tucson, AZ 85721, USA}
\IEEEauthorblockA{ \IEEEauthorrefmark{2} Applied Mathematics and Plasma Physics,\\ Los Alamos National Laboratory, Los Alamos, NM 87519, USA}
\IEEEauthorblockA{\IEEEauthorrefmark{3} NOGA, Power System Operator of Israel, Haifa, Israel}
}

\begin{document}
\maketitle
\begin{abstract}
Static contingency screening does not identify all of the transient
exposure that may follow a common transmission fault, particularly when
the affected element and the overloaded element are different. We
examine this problem using a linear stochastic swing surrogate with
short-lived, fault-local multiplicative variability. A cumulative
overload indicator distinguishes persistent exposure from brief
crossings, and categorical--exponential importance sampling estimates
its exceedance probabilities under a fixed nominal scenario law.
On the IEEE 118-bus system, the stochastic and deterministic formulations
produce different patterns of predicted exposure, with more monitored
lines affected in the reported stochastic comparisons. The sampling
tables report lower CEM runtimes at their selected budgets, but do not
establish equal statistical precision or an end-to-end real-time update.
The results demonstrate sensitivity to the prescribed uncertainty
representation and identify contingencies and monitored elements for
further assessment. They do not independently establish that the
stochastic probabilities are physically more accurate than those from
the deterministic model.

\end{abstract}
\begin{IEEEkeywords}
Importance Sampling, Instantons for Rare Impactful Events, N-1 Contingencies, Power System Dynamics under Faults, Fluctuations, Transmission Grids.
\end{IEEEkeywords}
\section{Introduction}
\label{sec:intro}
\subsection{Motivation}
Power-system stability is concerned with the ability of the grid to regain an acceptable operating state after disturbances without loss of synchronism or cascading failures. Of particular relevance is transient stability, which governs the short post-fault interval where electromechanical dynamics, protection actions, and topology changes interact most strongly \cite{machowski2020power}. Increasing heterogeneity and uncertainty in modern grids motivate real-time dynamic security assessment as a complement to traditional offline studies \cite{xue2006power}. Recent large-scale events demonstrate that fast fault-driven interactions can trigger low-probability, high-impact outcomes that static screening fails to capture \cite{entsoe2025iberian,fercnerc2021uri}.

From an operator's perspective, not all $N-1$ contingencies are equally informative. If a particular transmission line -- when faulted -- repeatedly appears among the highest-risk scenarios leading to transformer overcurrent, that line deserves increased operational scrutiny, even if it is not traditionally classified as critical. Likewise, assets -- and especially transformers -- that systematically accumulate post-fault stress become natural targets for monitoring and preventive control.

This work develops an online decision-support framework that quantifies, in probabilistic terms, the risk associated with counterfactual transmission faults. We focus on single-phase disturbances and on the likelihood of secondary effects during the seconds-long post-fault transient, such as overloads of critical equipment, e.g. transformers. Building on \cite{almada2025real}, which employs a deterministic linear swing model, we introduce fault-induced dynamic uncertainty by modeling short-lived, localized fluctuations in the parameters of the faulted line\footnote{Fault-induced (co-occurring) uncertainty refers to transient, localized variability during a fault, distinct from ambient system-wide noise.}. 

The approach is demonstrated on the IEEE 118-bus model of the American Electric Power (AEP) system, aiming at scalable probabilistic transient security assessment for large transmission networks\footnote{Not a replacement for detailed dynamic simulation, the framework enables fast transient-based screening.}. 
\subsection{Main Contributions}
\begin{figure}[thb]
    \centering
    \includegraphics[width=1.0\columnwidth]{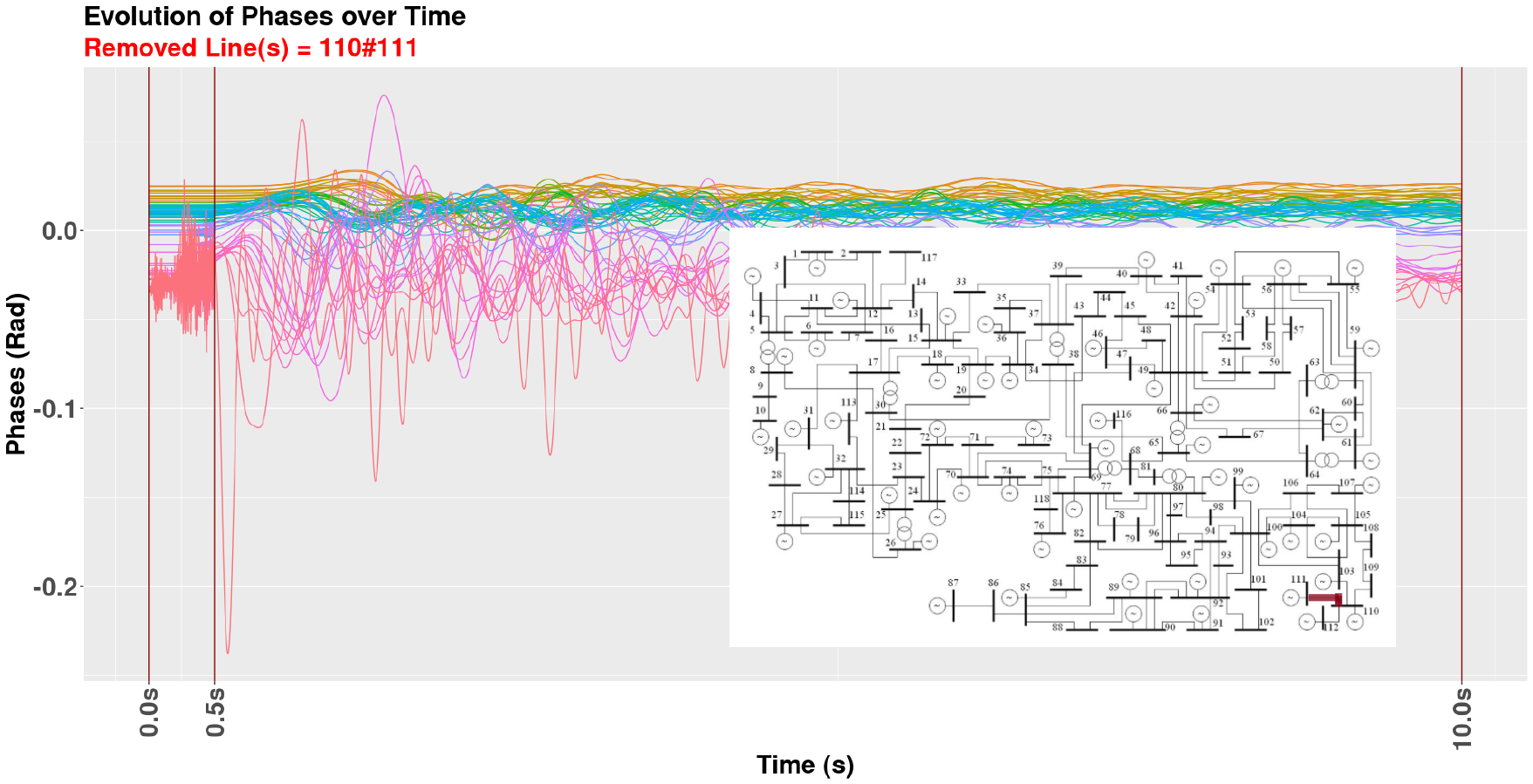}
	\caption{Example of a single-phase transmission fault whose post-clearing transient leads to severe downstream stress. While the faulted line is re-energized, secondary overloads on other elements may persist --motivating the need for dynamic $N-1$ screening beyond static criteria.} 
	\label{Point3}
\end{figure}

This work introduces the following contributions:
\begin{itemize}
\item \textbf{Efficient $N-1$ dynamic screening methodology} --- a scalable approach to analyze ensembles of transient responses associated with counterfactual faults on any transmission line, conditioned on a given operating state. Fig.~(\ref{Point3}) illustrates an example of undesirable transient behavior triggered by a $0.5\,\mathrm{s}$ line fault.

\item \textbf{Probabilistic risk framework} --- a stochastic modeling and estimation framework that characterizes the temporal evolution of system overloads driven by phase-angle dynamics, enabling quantitative assessment of grid reliability under fault-induced uncertainty.
\end{itemize}

\begin{tcolorbox}[
  colback=gray!5,
  colframe=black,
  boxrule=0.8pt,
  arc=2pt,
  left=6pt,
  right=6pt,
  top=6pt,
  bottom=6pt
]
\textbf{Operator View --- What this tool provides}

\vspace{2mm}

For a given operating state, the proposed dynamic N$-1$ screening addresses three questions relevant to an eventual operational screening tool:
\begin{itemize}
  \item \textbf{Which faulted transmission lines} most frequently trigger dangerous post-fault overloads elsewhere in the grid.
  \item \textbf{Which transformers or critical lines} repeatedly appear overloaded across the top-ranked fault scenarios.
  \item \textbf{How likely} these overloads are to occur during the first seconds following fault clearance.
\end{itemize}

This information enables operators to focus attention on specific lines or transformers during stressed conditions, even when those elements are not flagged as critical by conventional N$-1$ analysis.

\end{tcolorbox}

The framework separates reduced trajectory evaluation from targeted contingency sampling\footnote{A counterfactual contingency is a simulated event evaluated under the current operating state, not a fault currently active in the network.}. These components address different computational costs:
\begin{enumerate}
\item \textbf{Reduced dynamics:} linear stochastic propagation during the modified coupling interval, followed by deterministic matrix-exponential propagation after reclosure. 

\item \textbf{Targeted sampling:} use of cross-entropy--based importance sampling to focus computational effort on low-probability but highest-impact contingency scenarios.
\end{enumerate}
\subsection{Outline}
The remainder of this paper is organized as follows (Sections \textbf{\ref{sec:model}}--\textbf{\ref{sec:conclusion}}):
\begin{itemize}
\item \textbf{System and Fault Modeling} [Section \textbf{\ref{sec:model}}]: \\
We develop a scalable linearized second-order stochastic model of transmission-level dynamics based on the swing equations, explicitly incorporating co-occurring noise localized on the faulted line(s). Single-phase transmission faults are modeled over three regimes---pre-fault, fault-on, and post-clearing---capturing both deterministic and stochastic transient effects.
\item \textbf{Overload Indicator} [Section \textbf{\ref{sec:indicator}}]:\\
We introduce a continuous overload indicator that quantifies fault severity by combining thermal limits, fault duration, and phase-angle oscillations, enabling a time-resolved assessment of transient system stress.
\item \textbf{Scalable Dynamic Evaluation} [Section \textbf{\ref{sec:algorithm}}]:\\
We derive analytic and semi-analytic reductions of the governing dynamics that enable efficient re-evaluation of fault-induced transients across large ensembles of counterfactual contingencies.
\item \textbf{Probabilistic Risk Estimation} [Section \textbf{\ref{sec:mcmc}}]:\\
We estimate rare-event probabilities using Monte Carlo techniques and an efficient Cross-Entropy--based importance sampling method, compared with independent direct MC, and implemented within our custom simulator, \textbf{N1Plus}.
\item \textbf{Case Study: IEEE 118-Bus System} [Section \textbf{\ref{sec:casestudy}}]:\\
The proposed framework is evaluated on the IEEE 118-bus test system (118 buses, 177 lines, 19 generators, 9 transformers, and 91 loads), examining the consequences of the prescribed stochastic correction and the reported sampling costs.
\item \textbf{Conclusion and Outlook} [Section \textbf{\ref{sec:conclusion}}]:\\
We summarize the main findings and discuss extensions to larger and more heterogeneous grids, including the integration of data-driven and AI-assisted tools to further scale the methodology.
\end{itemize}

Appendix~\ref{app:localized_noise_foundation} provides the analytical
foundation for the stochastic model, including its diffusion structure,
exact solution representation, time-resolved sampling constructions,
and conditional moment equations.
The overall structure and computational pipeline of the proposed approach are illustrated in Fig.~(\ref{fig:pipeline}):
\begin{figure}[thb]
  \centering
  \includegraphics[width=1.0\columnwidth]{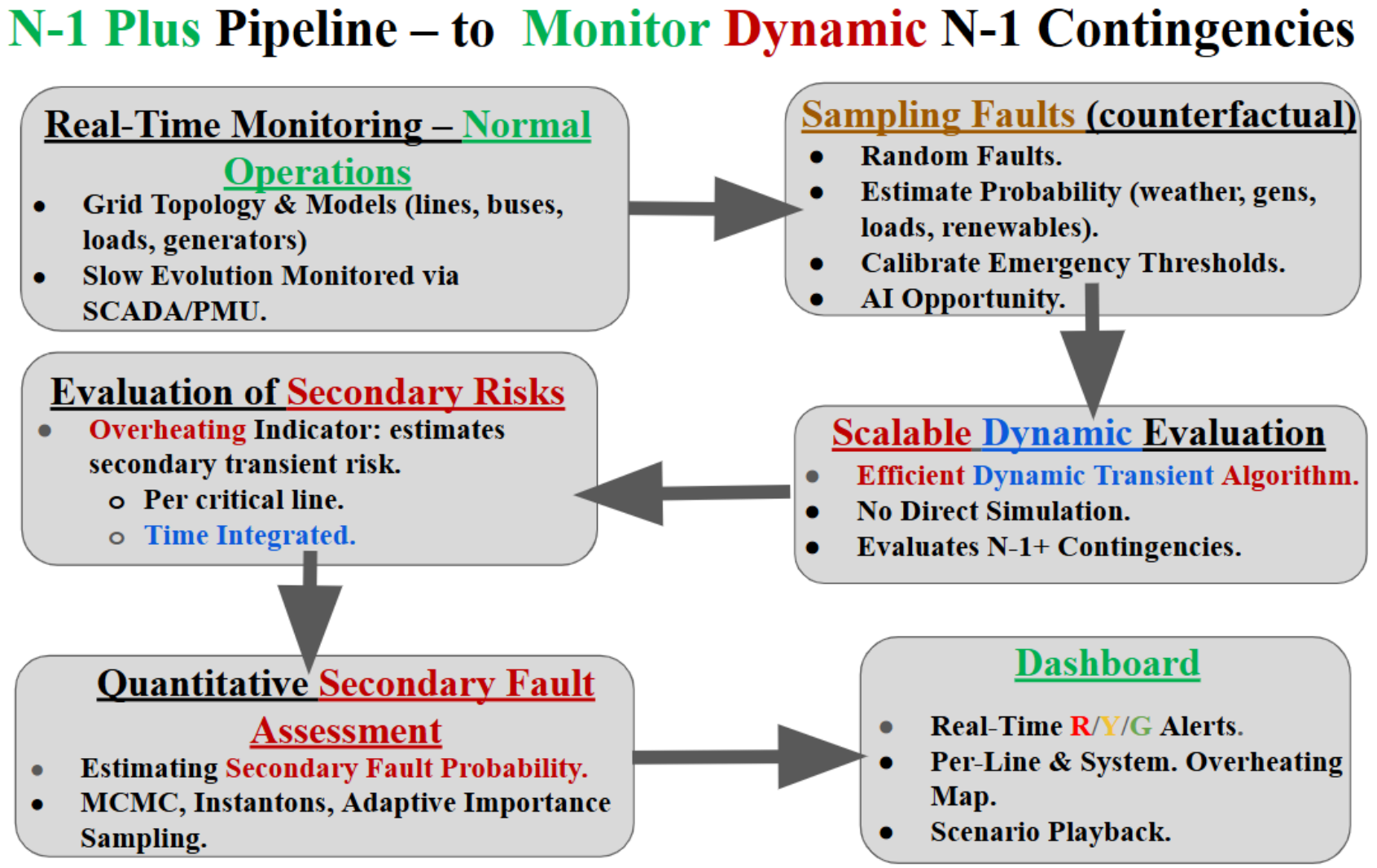}
  \caption{Conceptual Dashboard for Dynamic Contingency Screening.}
  \label{fig:pipeline}
\end{figure}
Appendix~\ref{app:p2_milstein_analysis} gives the algebraic form,
normalization, and discrete moments of the Milstein sampler.
Appendix~\ref{app:p2_experiment_record} records the supplied CEM
configuration, evaluation counts, and limits of the retained
experimental summaries.

\section{Related Work and Background}
\label{sec:related}
\subsection{$N-1$ Contingency Analysis}
$N-1$ contingency analysis\footnote{The proposed methodology also extends to broader contingency classes, including scenarios derived from historical fault data rather than strict $N-1$ rules.} is a cornerstone of power-system reliability assessment, requiring the grid to remain stable and within operational limits following the loss of any single critical component. This criterion underpins operational standards adopted by organizations such as NERC and ENTSO-E \cite{chatterjee2010n}. Conventional
$N-1$ studies rely on static power-flow simulations to identify limit violations and rank contingencies for corrective actions \cite{milano2010power}\footnote{Voltage regulation and dispatch corrections are not considered in this work.}. While effective in steady-state analysis, such approaches are computationally intensive and largely insensitive to transient dynamics, particularly in large, low-inertia systems with high renewable penetration. To address these limitations, we propose a framework for real-time screening of \emph{dynamic $N-1$ contingencies} that explicitly targets short-term post-fault behavior (sub-second to tens-of-seconds; see Table~\ref{tab:litreview}). Focusing on common disturbances such as single-phase short circuits\footnote{This paper focuses exclusively on single-phase transmission-line faults, which account for the majority of real-world events. Three-phase faults are intentionally excluded to maintain a tight operational scope and avoid conflating distinct protection and dynamic regimes.}, the framework estimates the probability of secondary effects -- specifically, transient but severe overloads on critical transmission elements -- that may precipitate further instability in stressed operating conditions.

In operational environments, contingency analysis results are typically presented as ranked lists of limit violations or security margins, rather than as explicit probability estimates. Control room tools such as EMS-based contingency analysis and real-time security assessment systems prioritize clarity, speed, and consistency with established reliability criteria, including NERC TPL standards and ENTSO-E operational practices \cite{chatterjee2010n,ENTSOE_Operational_Handbook,NERC_Reliability_Standards,wood2013power}. In this context, the probabilistic quantities computed in the proposed framework are \textbf{not intended for direct display}. Instead, they are used internally to \textbf{construct ranked risk categories and identify contingencies} that warrant further analysis. This representation is consistent with existing operational paradigms, where operators rely on severity indices, violation counts, and contingency rankings to guide decision-making. Accordingly, the proposed method should be interpreted as a screening layer that augments standard $N-1$ analysis by prioritizing contingencies based on their transient impact, while preserving compatibility with established tools and procedures. 
Despite extensive prior work on dynamic security assessment, most operational tools remain grounded in deterministic simulations of a limited set of contingencies, with emphasis on worst-case trajectories and security margins \cite{machowski2020power,xue2006power,milano2010power,ajjarapu2002bibliography,pavella2012transient,van2007voltage}. While these approaches are well established and widely deployed, they do not explicitly quantify the variability of post-fault responses or the likelihood of rare but high-impact transient events. In particular, conventional time-domain simulation evaluates individual trajectories and therefore provides limited information on how sensitive system behavior is to small perturbations in fault conditions, clearing times, or localized disturbances.
The present work differs in that \textbf{it considers an ensemble of fault realizations and evaluates the distribution of transient outcomes} using a stochastic formulation. This enables the identification of contingencies that are not necessarily critical under nominal conditions, but that exhibit \textbf{elevated risk due to variability in the post-fault response}. By combining analytic reductions of the swing equations with importance sampling techniques, the proposed approach provides a computationally efficient means of approximating this ensemble behavior, complementing existing deterministic methods rather than replacing them.
The goal is to provide a concise ranking of contingencies based on estimated transient overload risk, enabling prioritization within existing contingency analysis workflows.

In practice, the proposed method can be used upstream of standard $N-1$ analysis to identify a subset of contingencies for detailed time-domain simulation, thereby reducing computational burden while retaining compatibility with established reliability criteria.
\begin{table*}[!t]
\caption{N-1 Analysis vs Dynamic Evaluation.}
  \label{tab:litreview}
\small
\centering
\begin{tabular}{lll}
\hline
\multicolumn{1}{c}{\textbf{Aspect}}
 & \multicolumn{1}{c}{\textbf{N-1 Contingency}} & \multicolumn{1}{c}{\textbf{Swing Equations}} \vphantom{$\int_f^f$} \\
\hline
\vphantom{$\int_f^f$}Objective & Secure post-fault & Assess transient stability \vspace{1pt}\\
Time Scale [s] & $10^0 - 10^2$ & $10^{-3} - 10^{0}$ \vspace{3pt}\\
Application & Operations and planning& Dynamic stability  studies\vspace{1pt}\\
Model & AC/DC power flow &  Generator and load  responses as ODEs\vspace{1pt}\\
\hline
\end{tabular}
\end{table*}
\subsection{Rare Event/Instanton Sampling}
Rare--event (instanton) methods provide an efficient framework for identifying low-probability, high-impact events in complex systems \cite{juneja2006rare,asmussen2007stochastic,touchette2009large,rubino2009rare,freidlin2012random,hartmann2012efficient,grafke2015instanton}, including power grids \cite{dobson_complex_2007,chertkov2010predicting,chertkov_exact_2011,pfitzner_statistical_2011,owen2019importance}. Unlike standard Monte Carlo, these approaches leverage large-deviation theory to characterize the most probable rare trajectories -- \emph{instantons} -- that drive extreme outcomes such as blackouts and cascading failures \cite{dobson_complex_2007,pfitzner_statistical_2011}. The instanton minimizes a stochastic action (Freidlin--Wentzell) and thereby reveals both the physical mechanism and the exponential scaling of the associated probabilities \cite{juneja2006rare,freidlin2012random,grafke2015instanton}. Adaptive importance sampling (AIS), and in particular the cross-entropy (CE) method \cite{rubinstein_optimization_1997,de2005tutorial}, can be interpreted as using instanton-informed proposal families to bias sampling toward the most relevant regions of state space; iterative CE/AIS updates refine this bias via stochastic optimization, yielding substantial variance reduction and, in the low-probability regime, efficient estimators under suitable proposal families and sampling conditions \cite{rubinstein_optimization_1997,de2005tutorial,asmussen2007stochastic,hartmann2012efficient,owen_monte-carlo_2013}.
\section{Dynamics, Faults and Overloads}
\label{sec:model}
This section introduces the core modeling framework used throughout the paper: a transmission-level power grid dynamics model based on swing equations. We begin by formulating the system's behavior under small-to-moderate disturbances, then describe how line faults are incorporated into the model, and finally define a safety domain that encodes grid element overload constraints. This framework enables both theoretical analysis and efficient numerical simulation of dynamic contingencies.
\subsection{Modeling with Swing Equations}
We model the transmission network as a collection of generators, loads, and transmission lines operating on multiple time scales. The retained variables describe a reduced electromechanical response; the equations do not resolve independent device-level load, voltage-control, or protection dynamics. These omissions define the scope of the surrogate rather than a claim that physical load dynamics are always slow.

Each bus $i\in\mathcal{V}$ is associated with a complex voltage $v_i=V_i e^{j\theta_i}$, where $V_i$ denotes the voltage magnitude and $\theta_i$ the phase angle. Although transmission systems are inherently three-phase, we adopt the standard balanced-network assumption and represent each line by a single equivalent phase \cite{bergen2009power}. The grid is modeled as an undirected graph $\mathcal{G}=(\mathcal{V},\mathcal{E})$, where $\mathcal{V}$ denotes the set of buses and $\mathcal{E}$ the set of transmission lines. The reduced model uses a lossless network with fixed voltage magnitudes. This assumption neglects voltage magnitude variations and unbalanced conditions during faults, which can affect transient power flows and overload levels. More detailed models capture these effects through coupled voltage and dynamic equations \cite{kundur2007power}. As a result, the present formulation may misestimate transient stress in some scenarios, reflecting a trade-off between physical fidelity and computational efficiency appropriate for fast contingency screening.

Following \cite{almada2025real}, system dynamics are described by the linearized swing equations
\begin{align}
\label{eq:swing-lin}
m_i \ddot{\theta}_i + d_i \dot{\theta}_i 
+ \sum_{\{i,j\}\in\mathcal{E}} \beta_{ij}(\theta_i-\theta_j)=p_i,\qquad i\in\mathcal{V},
\end{align}
where $m_i$ and $d_i$ denote inertia and damping coefficients, $p_i$ represents net power injection, and $\beta_{ij}=V_iV_jB_{ij}$ is the effective stiffness of line $\{i,j\}$. In vector form,
\begin{align}
\label{irtrp}
\bm M \ddot{\bm\theta} + \bm D \dot{\bm\theta} + \bm L \bm\theta = \bm p,
\end{align}
with diagonal inertia and damping matrices $\bm M$ and $\bm D$, and weighted graph Laplacian $\bm L=\big(L_{ij}\big)$ defined by
\[
L_{ij}=\begin{cases}
\sum_{k\in\mathcal{V}}\beta_{ik}, & i=j,\\
-\beta_{ij}, & i\neq j.
\end{cases}
\]
Eq.~\eqref{irtrp} should be interpreted as a linearization of more general seconds-scale dynamics around a pre-disturbance operating point. During fault events, the parameters $\beta_{ij}$ associated with the affected line may change partially or completely, altering the network Laplacian and driving transient behavior. The modeling of these fault-induced variations is addressed in the following subsection.
\subsection{Effect of Co-Occurring Noise}\label{sec:nambientnoise}
The reduced model does not resolve the fast electrical and device
processes associated with a fault. We represent part of their possible
influence by a short-lived, line-local stochastic correction
\cite{wood2013power,sauer2017power}. The correction is prescribed in
this study rather than identified from EMT observations. Its purpose
is to examine sensitivity of the predicted electromechanical response
to localized variability, not to derive a white-noise law from arcing or
protection physics. The formal coupling perturbation in
\eqref{susuncert} motivates the location of the correction; the actual
phase-block SDE is specified in Sec.~\ref{sec:algorithm}.

\subsection{Simulating Faults}\label{sec:faults}
We consider faults occurring exclusively on transmission lines; faults at buses (e.g., generators or transformers) are excluded, although the framework readily extends to these cases. The assessment is conditional on a fault occurring; it does not estimate calendar-time fault-arrival rates. The nominal experiment selects an initiating contingency uniformly from the specified fault catalogue. In this work, we focus on \textbf{single-phase-to-ground faults}, which constitute the dominant class of transmission disturbances \cite{horowitz2022power}.

Under standard operating rules, the remaining two phases may continue to operate during a single-phase fault until clearance, typically within one second. Adopting a balanced-network approximation, a single-phase fault of duration $\tau$ is modeled by a temporary reduction of the affected line's susceptance from $\beta_{ij}$ to $\tfrac{2}{3}\beta_{ij}$, reflecting reduced transfer capability; upon clearance, the susceptance is restored to its nominal value\footnote{This approximation neglects millisecond-to-second unbalanced dynamics; more detailed fault representations are left for future work.\label{footnote:disclaimer}}. The nominal duration $\tau$ is exponentially distributed with mean $\lambda_0=0.1\,\mathrm{s}$, equivalently rate $10\,\mathrm{s}^{-1}$, independently of the initiating line. This law describes the sampled modified-coupling interval; it is not a fitted distribution of relay settings.

\begin{figure}[thb]
  \centering
  \includegraphics[width=1.0\columnwidth]{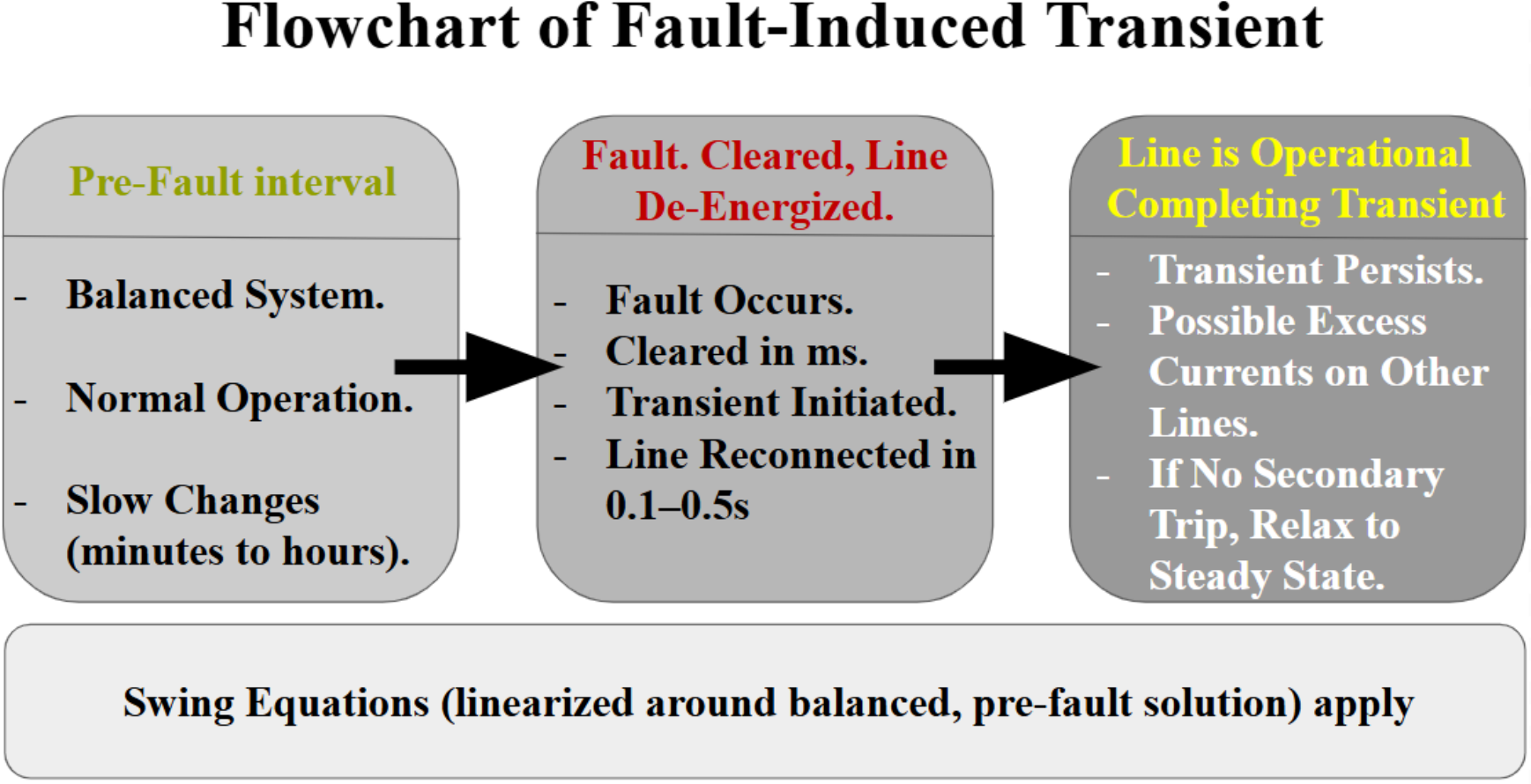}
  \caption{Fault timeline: system is balanced, pre-fault; fault is cleared, a faulty line is de-energized; line is back in service, post-fault transient.}
  \label{fig:timeline}
\end{figure}

Conditioned on a fault of duration $\tau$, a formal line-local representation of the coupling variability is
\begin{align}\label{susuncert}
\beta_{ij}(t)=\overline{\beta_{ij}}(t)+a_{ij}\zeta(t)\mathbb{I}\big(t\in[0,\tau]\big),
\end{align}
where $\zeta(t)$ is formal unit white noise, not a pointwise electrical measurement, $\overline{\beta_{ij}}(t)=\begin{cases}\tfrac{2}{3}\beta_{ij}, & t\le \tau,\\\beta_{ij}, & t>\tau,\end{cases}$ and $a_{ij}$ is a formal susceptance-noise amplitude, with units
$[\beta_{ij}]\,\mathrm{s}^{1/2}$. It is not the phase-diffusion
parameter $\sigma_{ij}$ in \eqref{eq:swing-lin-matrix}; no such
susceptance process is sampled in the reported runs.
The reported assignment $\sigma_{ij}=\beta_{ij}$ is a numerical scaling in the fixed unit convention used by the reduced model, not an equality of physical dimensions between a white-noise amplitude and a susceptance. It is not an EMT-fitted coefficient. This formulation captures the three stages of system evolution---\textcolor{mustard}{Pre-Fault}, \textcolor{red}{Fault}, and \textcolor{YellowOrange}{Operational} intervals---summarized schematically in Fig.~(\ref{fig:timeline}).

\subsection{Relation to Standard Fault Models}\label{sec:relstdfault}

The fault representation adopted in this work is intentionally simplified relative to detailed electromagnetic transient (EMT) and three-phase dynamic models commonly used in protection and stability studies. In standard transient stability analysis, faults are typically modeled as balanced three-phase short circuits followed by instantaneous line tripping, resulting in a temporary removal of the affected element from the network admittance matrix \cite{machowski2020power,horowitz2022power}. In contrast, real-world transmission faults are frequently single-phase-to-ground events, during which partial power transfer may persist until protection clears the fault, and system variables exhibit fast, localized fluctuations due to arcing phenomena and switching actions.

The model used in this paper should therefore be interpreted as a \textbf{reduced-order surrogate} designed for rapid screening of large ensembles of contingencies, rather than as a replacement for high-fidelity EMT simulation. The temporary reduction of line susceptance to a fraction of its nominal value during the fault interval provides a tractable approximation of diminished transfer capability, while preserving network connectivity and permitting a linear reduced description of the resulting dynamics. This approach is consistent with simplified representations used in classical transient stability studies, where faults are modeled through modifications of network parameters within the swing-equation framework \cite{chertkov_exact_2011}.

The stochastic term is a prescribed effective correction localized
at the faulted line's endpoints. Arcing, grounding, and protection
interactions motivate attention to that interval, but this study does
not derive a white-noise limit for those processes. In particular, the
phase-block diffusion in \eqref{eq:swing-lin-matrix} is not obtained by
substituting a random susceptance into the acceleration equation. The
comparison tests the consequences of this specified correction within
the reduced model.

We emphasize that the proposed formulation does not attempt to replicate detailed unbalanced or voltage-dependent dynamics. Instead, it provides a computationally efficient approximation of post-fault electromechanical behavior that is suitable for \textbf{screening, ranking, and prioritization of contingencies} within real-time operational workflows. More detailed validation against high-fidelity dynamic simulations and extension to voltage-coupled models remain important directions for future work.

\subsection{Safety Polytope}
We define the \emph{safety polytope} as the domain of phase configurations that ensure all line flows remain within their respective thermal (or angular) limits. Specifically, it characterizes the set of phase angles for which the power transmitted along every line does not exceed its safety threshold:
$$
\Pi_{\bm{\theta}} \!\doteq\! \left\{\! \bm{\theta} \!\in\! [-\pi,\pi]^{|\mathcal{V}|} \middle| \forall \{i,j\} \!\in\! \mathcal{E}:\left|\beta_{ij}\left(\theta_i \!-\! \theta_j\right)\right| \!\leq\! \bar{p}_{ij} \!\right\},
$$
where $\bar p_{ij}$ is the allowable power-flow magnitude in the
linearized relation, not an angle threshold. The set $\Pi_{\bm\theta}$
describes instantaneous compliance with those reduced flow limits. It
is not a certificate of transient stability, thermal safety, or protection
performance outside the assumptions of the linearized model.
\section{Overload Indicator}\label{sec:indicator}
To quantify \textbf{transient line stress}, we define an \textit{overload indicator} based on the solution of the linearized stochastic swing eqs.~\eqref{eq:swing-lin-matrix}. Let $\theta(t)=(\theta_i(t)\,|\,i\in\mathcal{V})$ denote the phase angles over $t\in[0,T]$, and let $\beta_{ij}(t)$ in this observable denote the deterministic effective coupling used for scoring: the prescribed reduced coupling during the open-phase interval and the nominal coupling afterward. The stochastic phase correction enters through \eqref{eq:swing-lin-matrix}; no formal white-noise value is inserted into a flow-threshold test. For a monitored line $\{i,j\}\in\mathcal{E}$, the \textbf{line-specific overload indicator} is
\begin{gather}
\label{eq:over-heat-ij}
S_{ij}(\theta_{0\to T}) \doteq \int_0^T 
\mathbb{I}\!\left(\left|\beta_{ij}(t)\big(\theta_i(t)-\theta_j(t)\big)\right|>\bar p_{ij}\right)\,dt,
\end{gather}
where $\bar p_{ij}$ is a safety threshold associated with the line's thermal limit and $\mathbb{I}(\cdot)$ is the indicator function. The monitored line need not coincide with the faulted line; overloads may be induced anywhere in the network by a given fault. To assess system-wide impact, we also define the \textbf{global
overload indicator} across a set of monitored lines $\mathcal{E}_{\rm m}\subseteq\mathcal{E}$:
\begin{gather}\label{glonOHI}
S(\theta_{0\to T}) \doteq \sum_{\{i,j\}\in\mathcal{E}_{\rm m}} S_{ij}(\theta_{0\to T}),
\end{gather}
with $\theta_{0\to T}\doteq\left(\theta_i(t)\,|\,i\in\mathcal{V},\,t\in[0,T]\right)$. The quantity $S(\theta_{0\to T})$ sums overload duration across monitored lines and can exceed $T$ when several lines are exposed. Its cumulative continuous-time integral is nondecreasing and continuous; numerical accumulation proceeds on the chosen observation grid. It records exposure duration, not peak flow or thermal damage.

The indicators are evaluated at the specified pre-fault operating state.
The probability law varies the initiating contingency and duration, together with the Wiener path.
The numerical score is denoted by $S_{h,\Delta}$, where $h$ describes
trajectory approximation and $\Delta$ the observation/quadrature grid.
The algorithm writes $S_k$ for this computed score. Accordingly,
$Q_{h,\Delta}(\gamma)=\mathbb P[S_{h,\Delta}>\gamma]$ is the target
of the implemented estimator; agreement of sampling methods does not
remove the difference between $Q_{h,\Delta}$ and $Q$.
Appendix~\ref{app:p2_duration_score} states the scoring and conditioning distinctions.

\section{Scalable Dynamic Evaluation}
\label{sec:algorithm}
Let $\bm x=(\bm\omega,\bm\theta)^{\mathsf T}$ contain the swing
frequency-deviation variables and phase coordinates. The deterministic
drift is obtained from \eqref{irtrp}; the fault-local phase correction is
specified separately. On the observation window,
\begin{align}\label{eq:swing-lin-matrix}
 d\bm x&=(\overline{\bm A}(t)\bm x+\bm P)\,dt
 +\mathbb I\{0\leq t\leq\tau\}G_{ij}\bm x\,dW_t,\\
 G_{ij}&=\begin{bmatrix}
 \mathbb O_{n\times n}&\mathbb O_{n\times n}\\
 \mathbb O_{n\times n}&\sigma_{ij}\bm M^{-1}v_{ij}v_{ij}^{\mathsf T}
 \end{bmatrix},\quad v_{ij}=e_i-e_j,\nonumber\\
 \bm P&=\begin{pmatrix}\bm M^{-1}\bm p\\\bm0\end{pmatrix},\qquad
 \bm x_0=\begin{pmatrix}\bm0\\\bm\theta_0\end{pmatrix},\nonumber\\
 \overline{\bm A}&(t)=\begin{bmatrix}
 -\bm M^{-1}\bm D&-\bm M^{-1}\overline{\bm L}(t)\\
 \mathds{1}_{n\times n}&\mathbb O_{n\times n}
 \end{bmatrix}\nonumber\\
 &=\bm A_0+\mathbb I\{0\leq t\leq\tau\}\,\delta\bm A.\nonumber
\end{align}
The initial phase vector solves $\bm L_0\bm\theta_0=\bm p$ with a fixed
additive phase reference. This does not require an inverse of the
ungrounded Laplacian. The inverse of $\bm M$ applies only to an
invertible retained differential model; it is not a prescription for
inverting zero-inertia physical load rows.

Equation~\eqref{eq:swing-lin-matrix} is interpreted in the It\^o sense.
Its lower-right diffusion block specifies the effective phase correction. Consequently, during the stochastic
interval $d\bm\theta=\bm\omega\,dt+
\sigma_{ij}\bm M^{-1}v_{ij}v_{ij}^{\mathsf T}\bm\theta\,dW_t$:
$\bm\omega$ is not a pathwise derivative of the noisy phase. When the
correction is absent, $\bm\omega=\dot{\bm\theta}$ and the deterministic
swing equation is recovered. A literal white-noise perturbation of the
Laplacian in \eqref{irtrp} would instead enter the acceleration block.
The phase correction here is therefore an effective closure motivated
by fault locality, not an algebraic identity obtained by substituting
\eqref{susuncert} into the second-order equation.

Let $\widetilde{\bm A}=\bm A_0+\delta\bm A$. The stochastic fundamental
matrix is defined by
\[
 d\Phi(t,s)=\widetilde{\bm A}\Phi(t,s)\,dt
            +G_{ij}\Phi(t,s)\,dW_t,\qquad \Phi(s,s)=I.
\]
Variation of constants gives
\begin{gather}\label{soldat}
\begin{aligned}
 \bm x(t)&=\Phi(t,0)\bm x_0+\int_0^t\Phi(t,s)\bm P\,ds,
 &&0\leq t\leq\tau,\\
 \bm x(t)&=e^{\bm A_0(t-\tau)}\bm x_\tau
 +\int_\tau^t e^{\bm A_0(t-s)}\bm P\,ds,
 &&t>\tau.
\end{aligned}
\end{gather}
This representation is exact for the stated SDE, but it is not by itself
an algorithm for evaluating $\Phi$. In the special commuting case
$[\widetilde{\bm A},G_{ij}]=0$, it reduces to
\[
 \Phi(t,s)=\exp\!\left[
 (\widetilde{\bm A}-\tfrac12G_{ij}^2)(t-s)
       +G_{ij}(W_t-W_s)\right].
\]
Without that commutation property, the single-exponential expression is
not an exact propagator; the general stochastic flow must be evaluated
with a method that accounts for the noncommuting operators
\cite{r4a,r4b}. No commutation assumption is imposed on the network
operators in this paper.

For the diffusion specified in \eqref{eq:swing-lin-matrix},
\eqref{eq:noiseapp:commutator} shows that the drift and diffusion do not
commute whenever $\sigma_{ij}\ne0$.

A change to eigenvector coordinates also retains the complete diffusion
transformation. If
$\widetilde{\bm A}=\widetilde{\bm U}\widetilde{\bm\Lambda}
\widetilde{\bm U}^{-1}$ and $\bm A_0=\bm U\bm\Lambda\bm U^{-1}$,
then
\[
\begin{aligned}
 d\bm\xi_{\rm f}
 &=\left(\widetilde{\bm\Lambda}\bm\xi_{\rm f}
       +\widetilde{\bm U}^{-1}\bm P\right)dt+\widetilde{\bm U}^{-1}G_{ij}\widetilde{\bm U}
       \bm\xi_{\rm f}\,dW_t,&&t\leq\tau,\\
 \dot{\bm\xi}_0&=\bm\Lambda\bm\xi_0+\bm U^{-1}\bm P,
 &&t>\tau.
\end{aligned}
\]
The initial coordinates are
$\bm\xi_{\rm f}(0)=\widetilde{\bm U}^{-1}\bm x_0$ and
$\bm\xi_0(\tau)=\bm U^{-1}\bm x_\tau$.
Diagonalizing the drift does not generally diagonalize the diffusion or
make the stochastic modes independent. Spectral evaluation of the
post-reclosure deterministic interval can be reused across samples
\cite{pagnier2019optimal}; the stochastic propagation cost remains a
separate part of each sampled trajectory.

Appendix~\ref{app:noise_flow} derives \eqref{soldat} and its
integrating-factor representation. The stochastic trajectories reported
in this study are sampled with the scalar-driver Milstein method in
Appendix~\ref{app:noise_sampling}. Exponential splitting is retained
there as an analytical alternative, not as the method used to generate
the reported results. Milstein approximates the stochastic flow; the
variation-of-constants formula is not an exact finite-step sampling rule.
The method retains the fault-endpoint coupling and the random state
attained at reclosure. Appendix~\ref{app:p2_milstein_analysis} derives
its discrete moment recursions and distinguishes them from the exact
continuous moments.

The conditional mean and covariance satisfy the closed equations in
Appendix~\ref{app:noise_moments}. For the fixed It\^o drift, the mean
follows the deterministic baseline, while covariance and nonlinear
overload probabilities can depend on the noise amplitude. The moment
equations provide deterministic references for checking a numerical
sampler; they do not replace the temporally consistent trajectories
required by \eqref{eq:over-heat-ij}.

\section{Probabilistic Risk Assessment}
\label{sec:mcmc}
The overload indicators in \eqref{eq:over-heat-ij} and \eqref{glonOHI}
define the events to be estimated. For a threshold $\gamma\geq0$, let
\begin{equation}
\begin{aligned}
 Q_{ij}(\gamma)&=\mathbb P_{p_Z}\!\left[S_{ij}(\theta_{0\to T})>\gamma\right],\\
 Q(\gamma)&=\mathbb P_{p_Z}\!\left[S(\theta_{0\to T})>\gamma\right].
\end{aligned}
\label{eq:strict_risk_probabilities}
\end{equation}
The comparison is strict throughout: an exposure equal to the threshold
is not counted. In particular, $Q(0)$ is the probability of positive
accumulated overload, not the probability of the certain event $S\geq0$.
The line-specific indicator measures accumulated overload duration;
the global indicator sums that duration over monitored elements and
can exceed the observation horizon. Neither quantity is a thermal-damage
model or a probability of loss of synchronism.

\subsection{Cross-Entropy Method (CEM)}\label{subsec:cem}
\paragraph{Nominal contingency law.}
Let $\mathcal F$ be the finite catalogue of candidate initiating
contingencies, and let $E=|\mathcal F|$. An entry identifies the line
whose coupling is modified; the catalogue and its mapping to network
branches are part of the case specification. This distinguishes sampled
identifiers from the monitored set $\mathcal E_{\rm m}$ and from an
aggregate count of network elements. With $Z=(e,\tau)$, the nominal law is
\begin{equation}
 p_Z(e,\tau)=\frac{1}{E\lambda_0}\exp\!\left(-\frac{\tau}{\lambda_0}\right),
 \qquad e\in\mathcal F,\quad \tau\geq0,
 \label{eq:nominal_contingency_law}
\end{equation}
where $\lambda_0=0.1\,\mathrm{s}$ is the exponential mean. Consequently,
$\mathbb E[\tau]=0.1\,\mathrm{s}$ and the rate is
$10\,\mathrm{s}^{-1}$. Initiating-category selection and duration are
independent. These probabilities are conditional on an initiating fault
in $\mathcal F$; they do not include a no-fault probability or a
calendar-time occurrence rate. The duration law defines the nominal
scenario ensemble, not a fitted distribution of protection settings.

\paragraph{Proposal family and initialization.}
The proposal retains the line-indexed exponential means of the
original formulation:
\begin{equation}
\begin{aligned}
 q_\nu(e,\tau)&=\frac{\phi_e}{r_e}\exp\!\left(-\frac{\tau}{r_e}\right),\\
 \nu&=((\phi_e,r_e)_{e\in\mathcal F}),\\
 \phi_e&>0,\qquad \sum_{e\in\mathcal F}\phi_e=1,\qquad r_e>0.
\end{aligned}
\label{eq:categorical_exponential_proposal}
\end{equation}
The restriction $r_e\equiv r$ gives a shared-mean family. All duration
parameters are exponential means measured in seconds; $\lambda_0$
remains fixed. The supplied scalar $r^{(0)}$ specifies the common
initial value $r_e^{(0)}=r^{(0)}$. It does not by itself determine
whether the means remain tied during subsequent updates.
The recorded initialization and adaptation budget are
\begin{equation}
\boxed{
\begin{gathered}
 N_{\rm CE}=200,\qquad K=15,\qquad \rho=0.1,\\
 E=177,\qquad n_{\rm elite}=20,\qquad \kappa=0.9,\\
 \phi_e^{(0)}=\frac1E,\qquad r^{(0)}=1.0\,\mathrm{s}.
\end{gathered}}
\label{eq:cem_initialization}
\end{equation}
Here $E=177$ is the number of transmission lines in the initiating-fault
catalogue. Nominal line selection and proposal initialization assign
probability $1/177$ to each line. The monitored elements are specified
separately and do not increase this initiating-fault count.
Appendix~\ref{app:p2_experiment_record} records the sampling configuration
and evaluation counts.

\paragraph{Elite selection and updates.}
At iteration $\ell$, draw $N_{\rm CE}$ scenarios from
$q_{\nu^{(\ell-1)}}$ and compute their global overload scores.
Let $\mathcal I_\ell$ contain the $n_{\rm elite}=\rho N_{\rm CE}$
largest scores, with ties resolved to retain exactly that number.
The proposal is fitted by maximizing the unweighted elite
log-likelihood \cite{rubinstein_optimization_1997,de2005tutorial}:
\begin{equation}
\begin{aligned}
 \widehat\phi_e^{(\ell)}
   &=\frac{1}{n_{\rm elite}}
      \sum_{k\in\mathcal I_\ell}\mathbb I\{e_k=e\},\\
 \widehat r_e^{(\ell)}
   &=\frac{\sum_{k\in\mathcal I_\ell}\tau_k\mathbb I\{e_k=e\}}
          {\sum_{k\in\mathcal I_\ell}\mathbb I\{e_k=e\}}.
\end{aligned}
\label{eq:elite_mean_updates}
\end{equation}
This duration update applies to categories represented in the elite
set; otherwise their previous means are retained. Under the tied
restriction, the update is instead the pooled mean
$\widehat r^{(\ell)}=n_{\rm elite}^{-1}\sum_{k\in\mathcal I_\ell}\tau_k$.
Neither update takes the reciprocal of the mean. Smoothing gives
\begin{equation}
\begin{aligned}
 \phi_e^{(\ell)}&=0.9\phi_e^{(\ell-1)}
                    +0.1\widehat\phi_e^{(\ell)},\\
 r_e^{(\ell)}&=0.9r_e^{(\ell-1)}+0.1\widehat r_e^{(\ell)}.
\end{aligned}
\label{eq:proposal_support_update}
\end{equation}
The coefficient $\kappa=0.9$ is the fraction of the previous value
retained. In exact arithmetic,
$\phi_e^{(\ell)}\geq0.9^{\ell}/E>0$; a category absent from an elite set
therefore keeps positive proposal probability. Adaptation stops after
exactly $K=15$ iterations. The elite cutoff guides proposal learning;
it is not the final exposure threshold $\gamma$, and reaching
$\gamma$ is not the stopping condition.

This unweighted elite fit is distinct from a nominal-density-weighted
cross-entropy objective. It is a proposal-learning rule, not a proof of
optimality for the final event. The complete likelihood ratio is used
in final estimation regardless of the adaptation rule.

\paragraph{Frozen-proposal estimation.}
After iteration $K$, freeze $\phi_e^\ast=\phi_e^{(K)}$ and
$r_e^\ast=r_e^{(K)}$, and draw $N$ fresh scenarios. Their weights are
\begin{equation}
 w_k=\frac{p_Z(e_k,\tau_k)}{q_{\nu^\ast}(e_k,\tau_k)}
 =\frac{r_{e_k}^\ast}{E\phi_{e_k}^\ast\lambda_0}
   \exp\!\left[\left(\frac1{r_{e_k}^\ast}-\frac1{\lambda_0}\right)\tau_k\right].
 \label{eq:joint_likelihood_ratio}
\end{equation}
Both category selection and duration are corrected. In the shared-mean
case, substitute $r_{e_k}^\ast=r^\ast$. The weights must use the means
that generated the final draws; the initializer alone does not identify
them. The estimators are
\begin{equation}
\begin{aligned}
 \widehat Q(\gamma)&=\frac1N\sum_{k=1}^N
                     w_k\mathbb I\{S_k>\gamma\},\\
 \widehat Q_{ij}(\gamma)&=\frac1N\sum_{k=1}^N
                     w_k\mathbb I\{S_{ij}^{(k)}>\gamma\}.
\end{aligned}
\label{eq:fresh_batch_estimators}
\end{equation}
The weights are not self-normalized. For independent final draws from
the frozen proposal with the stated support, conditioning on the
adaptation history gives the nominal expectation of the scored event
\cite{owen2019importance}. When $S_k$ is obtained by numerical
propagation and quadrature, this statement applies to that numerical
score, not automatically to the continuous-time functional. A
fault-conditioned estimate includes the additional factor
$\mathbb I\{e_k=e\}/\pi_e$, with $\pi_e=1/E$; proposal probabilities
are not conditional risk scores.

Conditional on $(e,\tau)$, the Milstein construction in
Appendix~\ref{app:noise_sampling} samples the Wiener path independently
of the scenario-selection law. The conditional Wiener distribution is
unchanged under nominal and proposal sampling, so its likelihood factors
cancel from \eqref{eq:joint_likelihood_ratio}.

\paragraph{Direct Monte Carlo.}
Direct MC draws $e_k$ uniformly from the same catalogue $\mathcal F$
and draws $\tau_k$ independently from the exponential law with mean
$\lambda_0=0.1\,\mathrm{s}$. It uses the same dynamics, scoring rule,
observation horizon, and strict comparison, without adaptation or
likelihood weighting:
\begin{equation}
 \widehat Q_{\rm MC}(\gamma)=\frac1N\sum_{k=1}^N\mathbb I\{S_k>\gamma\}.
 \label{eq:direct_mc_estimator}
\end{equation}
These are independent draws, not a Markov chain. If $\tau_k>T$,
reclosure lies outside the observation window. The original sampled
duration remains in the likelihood ratio; it is not replaced by $T$.

\subsection{N1Plus}\label{subsec:n1plus}
N1Plus combines a reduced trajectory model, the overload observables,
and contingency sampling. For every sampled initiating line and duration,
the switching times are retained, the response is evaluated through the
observation horizon, and the monitored-line scores are accumulated.
The method estimates overload exposure under the specified model rather
than certifying network security. It extends the static contingency
screening perspective in \cite{chertkov2010predicting} to transient
observables and uses importance weighting to preserve the nominal law
\cite{owen2019importance}. Algorithm~\ref{alg:n1plus} separates proposal
learning from the fresh estimation batch.

\begin{algorithm}[t]
\caption{N1Plus with a categorical--exponential proposal}
\label{alg:n1plus}
\begin{algorithmic}[1]
\Require Reduced model; identified fault catalogue $\mathcal F$;
 nominal mean $\lambda_0=0.1\,\mathrm{s}$; horizon $T$;
 threshold $\gamma\geq0$; final size $N$; recorded duration-parameter
 tying rule.
\State Set $N_{\rm CE}=200$, $K=15$, $n_{\rm elite}=20$,
       and $\kappa=0.9$.
\State Set $E=|\mathcal F|$, $\phi_e^{(0)}=1/E$, and
       $r_e^{(0)}=1.0\,\mathrm{s}$ for each $e$.
\State Check that the catalogue cardinality matches the case/run record.
\For{$\ell=1,\ldots,K$}
 \State Draw $N_{\rm CE}$ pairs $(e_k,\tau_k)$ from
        \eqref{eq:categorical_exponential_proposal}.
 \State Draw independent scalar Wiener increments for each path.
 \State Evaluate the fault-on SDE by Milstein and continue the attained
        state after reclosure; accumulate the scores $S_k$.
 \State Select exactly $n_{\rm elite}$ largest scores.
 \State Apply \eqref{eq:elite_mean_updates} with the recorded tying
        rule, then \eqref{eq:proposal_support_update}.
\EndFor
\State Freeze every $\phi_e^\ast$ and duration mean $r_e^\ast$.
\State Draw $N$ fresh scenarios from the frozen proposal.
\State Draw independent Wiener paths for the final scenarios and use the
       same Milstein and observation settings as in adaptation.
\State Evaluate $S_k$ and $S_{ij}^{(k)}$ through $T$.
\State Apply the complete weights \eqref{eq:joint_likelihood_ratio}.
\State Estimate strict exceedances using \eqref{eq:fresh_batch_estimators}.
\State Record $KN_{\rm CE}$ adaptation evaluations separately from $N$.
\Return $\widehat Q$ and $(\widehat Q_{ij})_{\{i,j\}\in\mathcal E_{\rm m}}$.
\end{algorithmic}
\end{algorithm}

\section{Case Study: IEEE 118-Bus System}\label{sec:casestudy}
\begin{figure*}[t]
    \centering
    \begin{subfigure}[t]{0.48\textwidth}
        \centering
        \includegraphics[width=1.0\columnwidth]{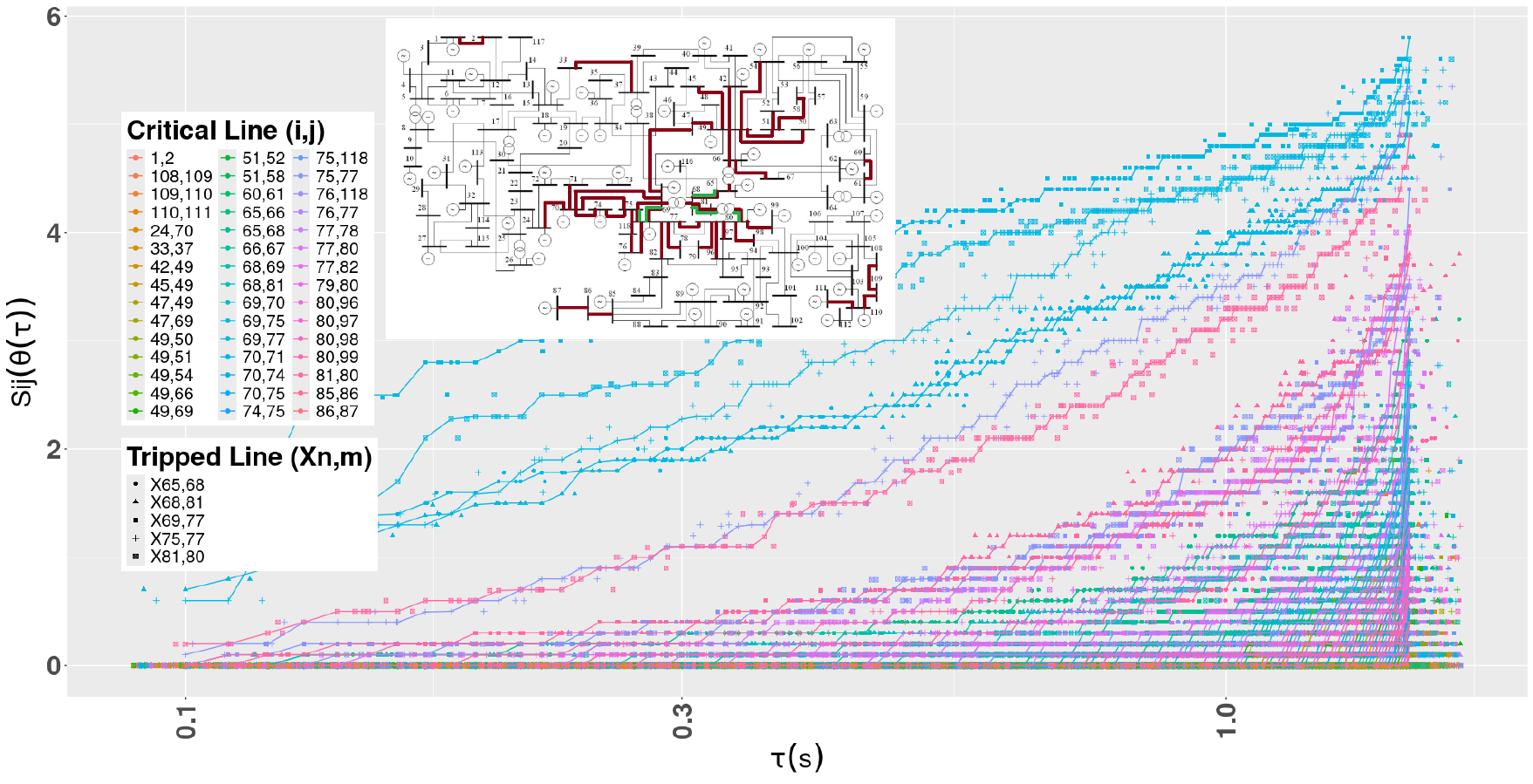}
        \caption{\textbf{Overload indicator vs. fault duration (single-phase faults)}. Each curve represents the overload of an individual line as a function of the fault duration $\tau$. \textbf{\textcolor{OliveGreen}{Tripped lines}} $(X_{n,m})$ identify initiating contingencies; \textbf{\textcolor{red}{critical transformer-connected lines}} $(i,j)$ are monitored elements with large accumulated overload duration. Relevant lines in the IEEE-118 network are highlighted.}
        \label{fig2}
    \end{subfigure}
    \hfill
    \begin{subfigure}[t]{0.48\textwidth}
        \centering
    	\includegraphics[width=1.0\columnwidth]{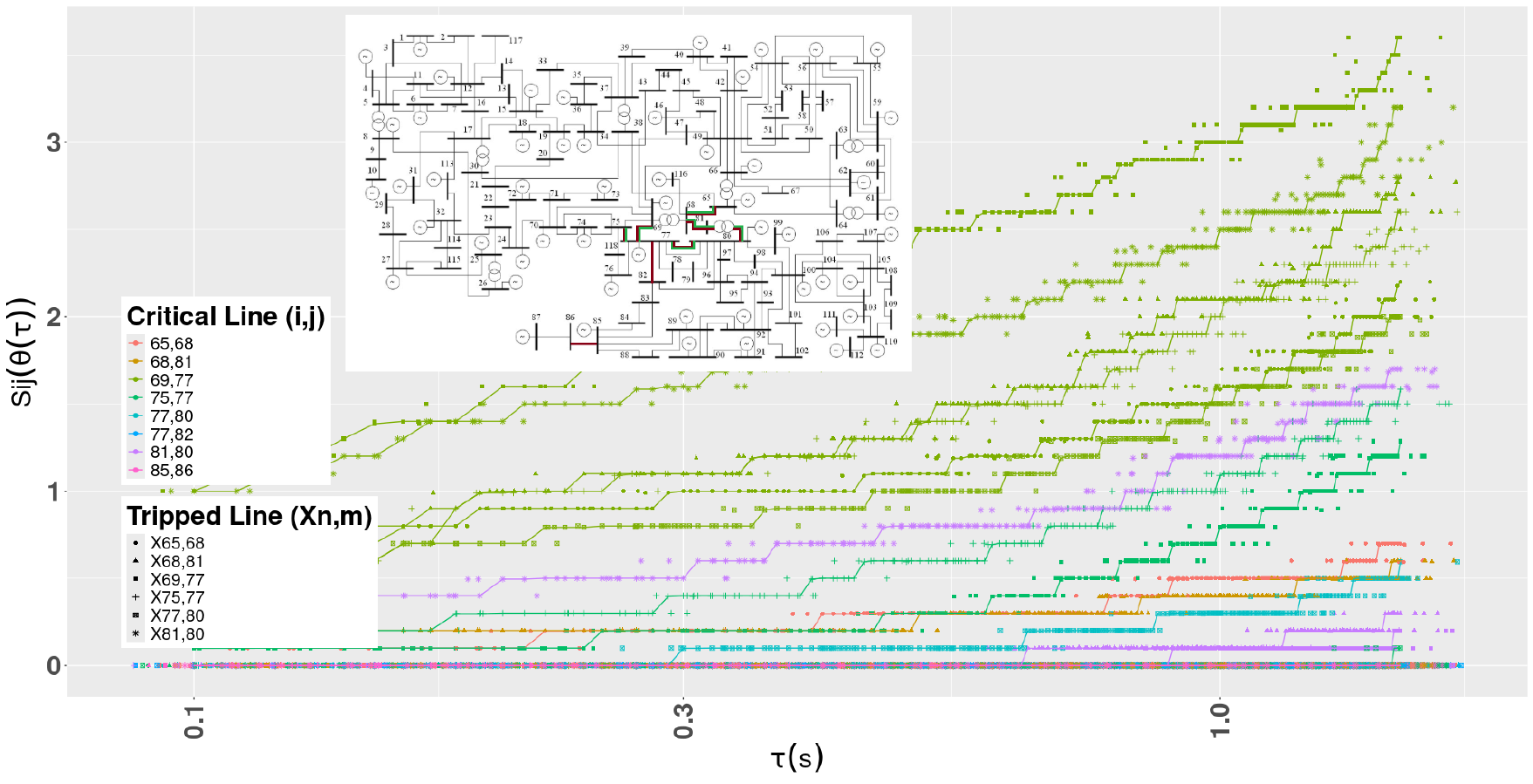}
    	\caption{\textbf{Overload indicator vs. fault duration without multiplicative noise (single-phase faults)}. Same definitions and notation as in Fig.~\ref{fig2}.}
    	\label{fig4}   
    \end{subfigure}
    \caption{\textbf{Overload evolution versus fault duration}. Comparison between stochastic dynamics (SDE, left) and deterministic dynamics (ODE, right).}
\end{figure*}

\begin{figure*}[t]
    \centering
    \begin{subfigure}[t]{0.48\textwidth}
        \centering
        \includegraphics[width=1.0\columnwidth]{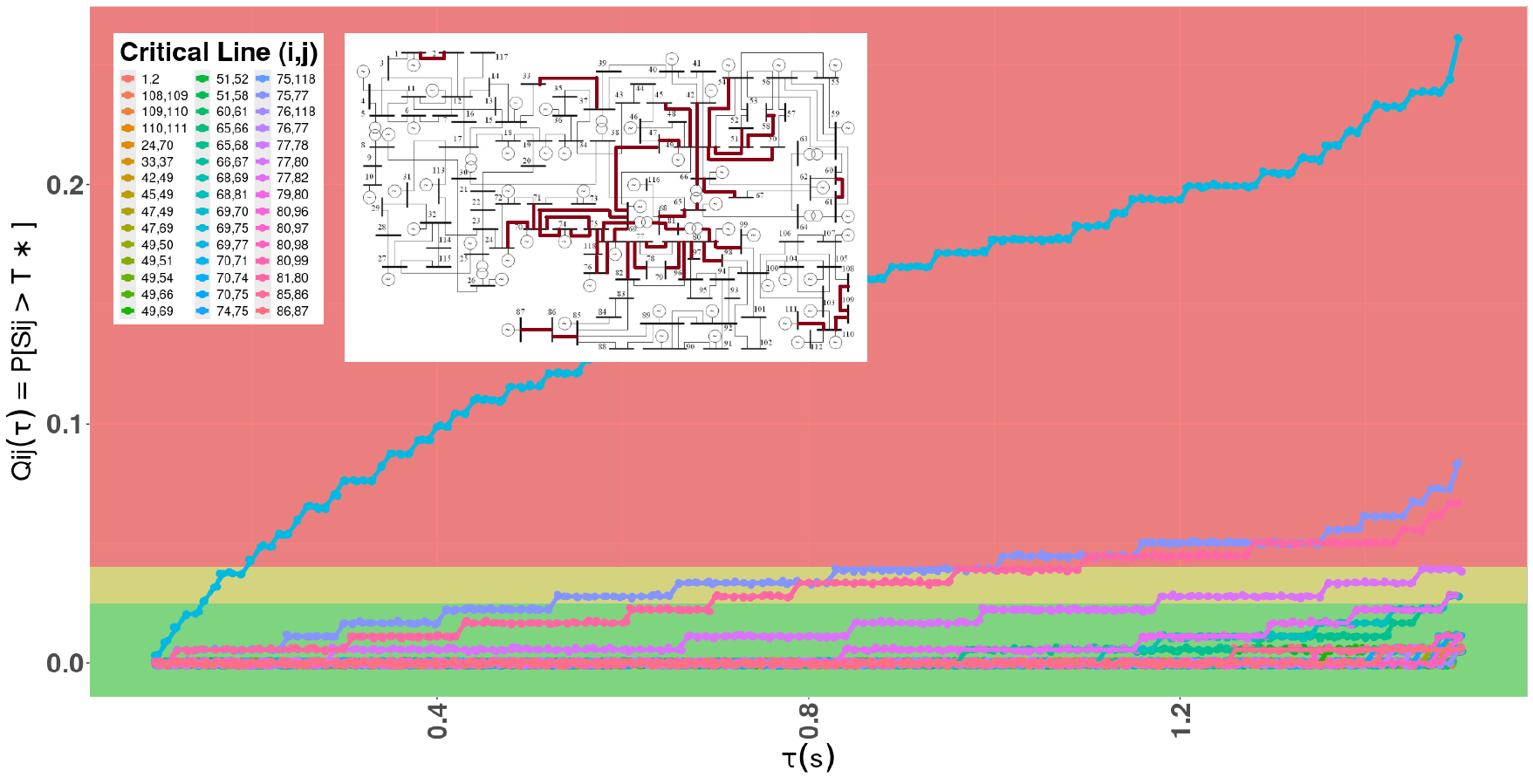}
        \caption{Retained duration-resolved probability-scale summaries for the stochastic model and strict threshold $T^{(*)}=1\,\mathrm{s}$. The bands are below $0.025$, from $0.025$ through $0.04$, and above $0.04$. The construction of the per-duration estimator is not identified in the retained plotting record; see Appendix~\ref{app:p2_duration_score}.}
        \label{fig3}
    \end{subfigure}
    \hfill
    \begin{subfigure}[t]{0.48\textwidth}
        \centering
    	\includegraphics[width=1.0\columnwidth]{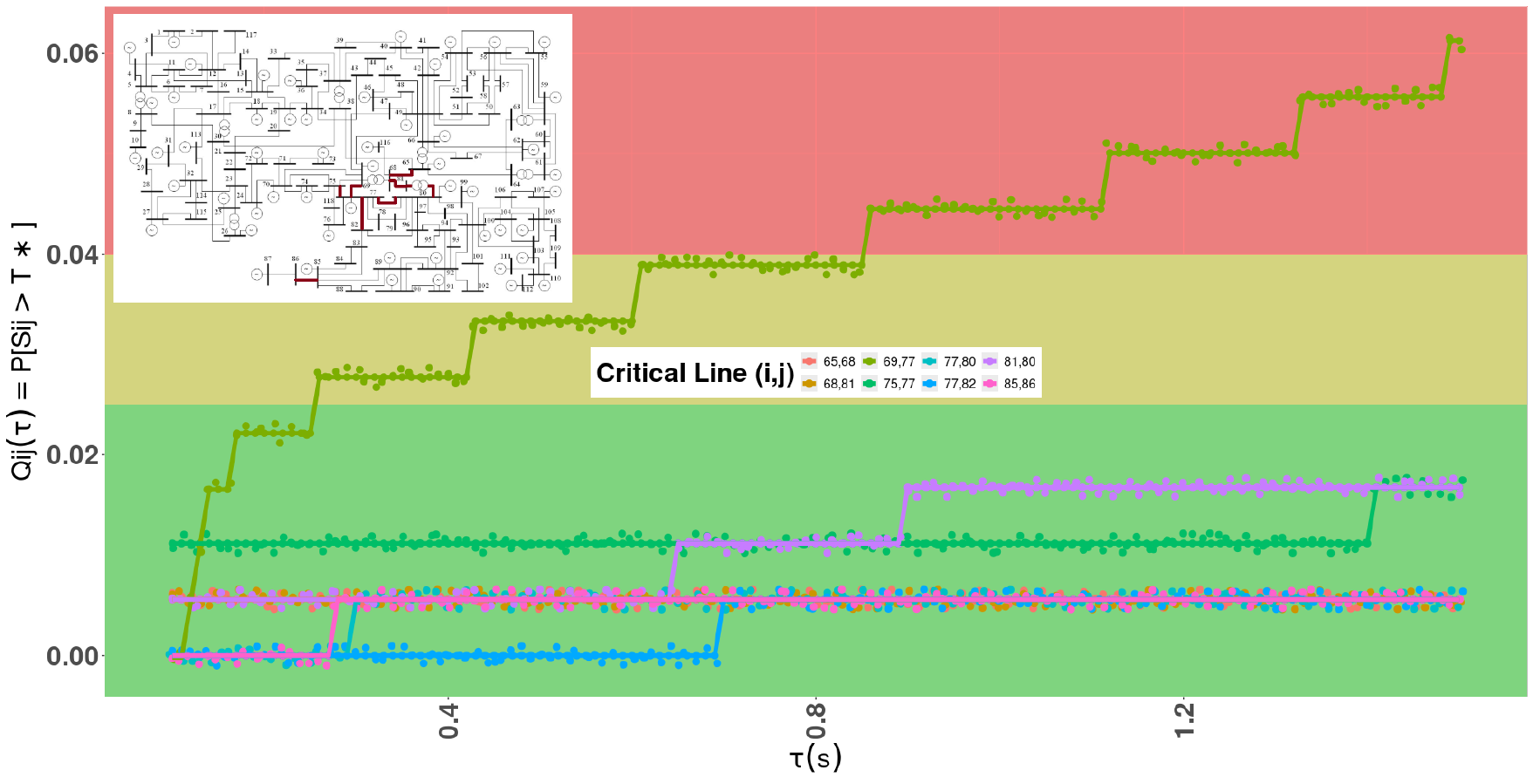}
    	\caption{Retained duration-resolved summaries for the deterministic comparator, with the same strict exposure threshold and illustrative bands as in Fig.~\ref{fig3}. The distinction between graphical scores and conditional probabilities is the same.}
    	\label{fig5}   
    \end{subfigure}
    \caption{Duration-resolved graphical summaries for the stochastic (left) and deterministic (right) descriptions. These plots are distinct from the duration-marginalized sampling tables.}
\end{figure*}
The reported case is the IEEE 118-bus system, with $19$ generators,
$35$ synchronous condensers, $177$ transmission lines, $9$ transformers,
and $91$ loads \cite{christie1993power}. These device counts do not
form a partition of buses: a bus may contain more than one type of
injection or device. The initiating-fault catalogue consists of the
$177$ transmission lines, with nominal probability $1/177$ per line.
The separate transformer count does not increase $E$; initiating faults
and monitored elements have distinct roles in the screening calculation.

The study compares the prescribed stochastic correction with its
zero-noise deterministic counterpart. Both use the stated single-phase
coupling change, exponential duration mean $0.1\,\mathrm{s}$, and
strict exposure threshold $T^{(*)}=1\,\mathrm{s}$. The stochastic runs
use Milstein during the noise-active interval. The comparison is at a
specified operating state; it does not include EMT-based parameter
calibration. Appendix~\ref{app:p2_experiment_record} distinguishes the
confirmed method and adaptation settings from the run details not
contained in the retained summaries. The deterministic comparator is

\begin{align}\label{eq:swing-lin-matrix2}
 \dot{\bm x}&=\bm A(t)\bm x+\bm P,\qquad
 \bm x=\begin{pmatrix}\bm\omega\\\bm\theta\end{pmatrix},\\
 \bm P&=\begin{pmatrix}\bm M^{-1}\bm p\\\bm0\end{pmatrix},\quad
 \bm x_0=\begin{pmatrix}\bm0\\\bm\theta_0\end{pmatrix},\quad
 \bm\omega=\dot{\bm\theta},\nonumber\\
 \bm A&(t)=\begin{bmatrix}
 -\bm M^{-1}\bm D&-\bm M^{-1}\bm L(t)\\
 \mathds{1}_{n\times n}&\mathbb O_{n\times n}
 \end{bmatrix}\nonumber\\
 &=\bm A_0+\mathbb I\{0\leq t\leq\tau\}\,\delta\bm A.\nonumber
\end{align}

The same threshold $T^{(*)}=1\,\mathrm{s}$ and risk classification are used for consistency. These thresholds are illustrative and are not intended to represent operational limits; in practice, such values should be specified by system operators based on empirical data and engineering judgment.
\subsection{Statistics of the Single Phase Fault with $T^{(*)}=1\text{s}$}\label{sec:3PF}
Figures~\ref{fig2} and~\ref{fig3} retain the reported stochastic
study, described by $N=100{,}000$ simulations. During each modified-
coupling interval, the affected line has two thirds of its nominal
coupling and the SDE is sampled by Milstein. The nominal duration law
has mean $0.1\,\mathrm{s}$; the per-duration construction of the figures
is a separate part of the plotting procedure.

The displayed monitored lines $\{69,77\}$, $\{81,80\}$, and
$\{75,77\}$ have recorded exposures above $1\,\mathrm{s}$ in the
longer-duration examples. These monitored elements need not be the
initiating lines. The corresponding probability-scale graphics cross
the illustrative $4\%$ band at longer displayed durations. Because the
per-duration estimator is not reconstructed, those crossings are
reported features of the graphic, not independently verified nominal
risk probabilities.

\subsection{Statistics of the Single Phase Fault with $T^{(*)}=1\text{s}$ and no noise in the system}\label{sec:1PF}
Figures~\ref{fig4} and~\ref{fig5} retain the deterministic comparison,
also described by $N=100{,}000$ simulations. The physical coupling
change, nominal duration law, and strict exposure threshold are the
same. The displayed $\{69,77\}$ response exceeds $1\,\mathrm{s}$ of
exposure for longer-duration examples. A common scenario law does not
by itself establish that the two experiments used paired random draws.

The reported selections contain $45$ monitored lines with larger
exposures in the stochastic setting and $8$ in the deterministic
comparison. The displayed extrema satisfy $S_{ij}>5.5\,\mathrm{s}$ in
Fig.~\ref{fig2} and $S_{ij}<3.5\,\mathrm{s}$ in Fig.~\ref{fig4}.
These counts describe the reported selections; a complete selection
rule and paired-scenario ledger are not available in the retained
record. They support a comparison of the displayed responses, not a
population-wide classification rate. The stochastic probability-scale
curves also reach larger displayed levels, but that observation is kept
separate from validation of their probability estimator.

\subsection{Benchmark Analysis}\label{BenchAn}

\begin{table*}[!t]
\caption{Reported sampling budgets, strict-exceedance summaries, and mean runtimes for the stochastic model. CEM and the nominal-sampling MC comparator target the law in \eqref{eq:nominal_contingency_law}. The entries do not establish matched precision or a convergence rate; see Sec.~\ref{BenchAn}.}
\label{table4}
\small
\centering
\begin{tabular}{rrrrrrr}
\hline
& \multicolumn{3}{c}{\textbf{Cross Entropy Method}}
& \multicolumn{3}{c}{\textbf{Monte Carlo Method}} \\
\multicolumn{1}{c}{$\boldsymbol{\gamma}$ \textbf{[s]}}
& \multicolumn{1}{c}{Sample Size}
& \multicolumn{1}{c}{$\boldsymbol{\hat{Q}\approx \mathbb{P}[S > \gamma]}$}
& \multicolumn{1}{c}{\textbf{Mean Time [s]}}
& \multicolumn{1}{c}{Sample Size}
& \multicolumn{1}{c}{$\boldsymbol{\hat{Q}\approx \mathbb{P}[S > \gamma]}$}
& \multicolumn{1}{c}{\textbf{Mean Time [s]}}
\vphantom{$\int_p^f$} \\
\hline
\vphantom{$\int_p^f$}
$\boldsymbol{0}$
& $2.000 \times 10^{3}$
& $\boldsymbol{2.426 \times 10^{-1}}$
& $\boldsymbol{1.862 \times 10^{2}}$
& $1.000 \times 10^{4}$
& $\boldsymbol{2.426 \times 10^{-1}}$
& $\boldsymbol{9.505 \times 10^{2}}$ \\
\vphantom{$\int_p^f$}
$\boldsymbol{0.5}$
& $3.500 \times 10^{3}$
& $\boldsymbol{1.314 \times 10^{-1}}$
& $\boldsymbol{3.598 \times 10^{2}}$
& $2.500 \times 10^{4}$
& $\boldsymbol{1.314 \times 10^{-1}}$
& $\boldsymbol{2.375 \times 10^{3}}$ \\
$\boldsymbol{5.0}$
& $1.000 \times 10^{4}$
& $\boldsymbol{1.311 \times 10^{-2}}$
& $\boldsymbol{9.479 \times 10^{2}}$
& $5.000 \times 10^{4}$
& $\boldsymbol{1.310 \times 10^{-2}}$
& $\boldsymbol{4.750 \times 10^{3}}$ \\
$\boldsymbol{10.0}$
& $1.250 \times 10^{4}$
& $\boldsymbol{2.121 \times 10^{-3}}$
& $\boldsymbol{1.235 \times 10^{3}}$
& $5.000 \times 10^{4}$
& $\boldsymbol{2.120 \times 10^{-3}}$
& $\boldsymbol{4.750 \times 10^{3}}$ \\
\hline
\end{tabular}
\end{table*}

\begin{table*}[!t]
\caption{Reported sampling budgets, strict-exceedance summaries, and mean runtimes for the deterministic comparator without multiplicative noise. The numerical entries are retained as reported. Their estimator and timing interpretation is qualified in Sec.~\ref{BenchAn}.}
\label{table5}
\small
\centering
\begin{tabular}{rrrrrrr}
\hline
& \multicolumn{3}{c}{\textbf{Cross Entropy Method}}
& \multicolumn{3}{c}{\textbf{Monte Carlo Method}} \\
\multicolumn{1}{c}{$\boldsymbol{\gamma}$ \textbf{[s]}}
& \multicolumn{1}{c}{Sample Size}
& \multicolumn{1}{c}{$\boldsymbol{\hat{Q}\approx \mathbb{P}[S > \gamma]}$}
& \multicolumn{1}{c}{\textbf{Mean Time [s]}}
& \multicolumn{1}{c}{Sample Size}
& \multicolumn{1}{c}{$\boldsymbol{\hat{Q}\approx \mathbb{P}[S > \gamma]}$}
& \multicolumn{1}{c}{\textbf{Mean Time [s]}}
\vphantom{$\int_p^f$} \\
\hline
\vphantom{$\int_p^f$}
$\boldsymbol{0}$
& $2.000 \times 10^{3}$
& $\boldsymbol{8.066 \times 10^{-2}}$
& $\boldsymbol{1.905 \times 10^{2}}$
& $1.000 \times 10^{4}$
& $\boldsymbol{8.066 \times 10^{-2}}$
& $\boldsymbol{9.429 \times 10^{2}}$ \\
\vphantom{$\int_p^f$}
$\boldsymbol{0.5}$
& $3.500 \times 10^{3}$
& $\boldsymbol{3.374 \times 10^{-2}}$
& $\boldsymbol{3.334 \times 10^{2}}$
& $2.500 \times 10^{4}$
& $\boldsymbol{3.374 \times 10^{-2}}$
& $\boldsymbol{2.431 \times 10^{3}}$ \\
$\boldsymbol{5.0}$
& $1.000 \times 10^{4}$
& $\boldsymbol{1.020 \times 10^{-5}}$
& $\boldsymbol{9.505 \times 10^{2}}$
& $5.000 \times 10^{4}$
& $\boldsymbol{9.850 \times 10^{-6}}$
& $\boldsymbol{4.829 \times 10^{3}}$ \\
$\boldsymbol{10.0}$
& $1.250 \times 10^{4}$
& $\boldsymbol{2.100 \times 10^{-6}}$
& $\boldsymbol{1.188 \times 10^{3}}$
& $5.000 \times 10^{4}$
& $\boldsymbol{1.890 \times 10^{-6}}$
& $\boldsymbol{4.891 \times 10^{3}}$ \\
\hline
\end{tabular}
\end{table*}

Tables~\ref{table4} and~\ref{table5} report the stochastic and
deterministic results from \eqref{eq:swing-lin-matrix} and
\eqref{eq:swing-lin-matrix2}, respectively. Both target
$Q(\gamma)=\mathbb P[S>\gamma]$ under the same exponential duration
law with mean $0.1\,\mathrm{s}$. The tested global thresholds are
$\gamma\in\{0,0.5,5,10\}\,\mathrm{s}$; their numerical values are
unchanged. The $\gamma=0$ row counts positive accumulated exposure.

Using the displayed mean runtimes, the MC-to-CEM ratios in
Table~\ref{table4} are approximately $5.105$, $6.601$, $5.011$,
and $3.846$. The corresponding ratios in Table~\ref{table5} are
$4.950$, $7.292$, $5.080$, and $4.117$. The two tables contain
different recorded timing values. Neither table separates adaptation
from final estimation or identifies the repetition count and timing
boundaries. The comparisons therefore describe the listed budgets;
they do not establish equal precision, a convergence rate, or a
controlled speedup relative to EMT \cite{kroese2012monte}.

The supplied adaptation budget contributes
$200\times15=3{,}000$ trajectory evaluations before a final batch.
Appendix~\ref{app:p2_experiment_record} separates this count from final
sample sizes. No adaptation wall time is inferred from a per-trajectory
average or added to the recorded times without a stage-specific log.

The MC probability column preserves the reported summaries rather than
reconstructing unrecorded event counts. In Table~\ref{table5}, the entries
$9.85\times10^{-6}$ and $1.89\times10^{-6}$ at $N=50{,}000$ are both
below $1/N=2\times10^{-5}$. They cannot be single-batch empirical
Bernoulli frequencies of the form \eqref{eq:direct_mc_estimator} at that
batch size. Rounding alone is not evidence of an estimator error, but it cannot
explain these positive values below the smallest positive single-batch
frequency. The table does not specify any repeated-batch
aggregation or other estimator construction that would account for this.
The values are retained as reported, without assigning single-batch
binomial uncertainty or using those entries to establish resolved
rare-event precision.

\subsection{Interpretation and Limits of the Comparison}\label{sec:qualvalcons}
The results in Figs.~\ref{fig2}--\ref{fig5} demonstrate sensitivity of
the modeled exposure to the prescribed stochastic correction. In the
reported comparisons, more lines accumulate substantial exposure and
the displayed probability-scale summaries reach higher levels than
those in the deterministic comparison. This is an observation about the two reduced descriptions,
not a general theorem that adding noise increases every overload
probability or that a deterministic calculation always understates
physical risk.

The overload indicator retains a linearized-flow interpretation only
within the assumptions of the underlying model. Voltage changes,
unbalanced electrical dynamics, protection actions, and large angle
separations can affect the physical consequences of the same initiating
fault \cite{machowski2020power,sauer2017power}. The present comparison
does not resolve those mechanisms. Agreement between MC and CEM when
both use the same surrogate concerns sampling, not validation of that
surrogate against a higher-fidelity reference.

\section{Conclusion and Outlook}
\label{sec:conclusion}
This study extends random-contingency screening by prescribing a
fault-local stochastic correction within the reduced dynamics. The
IEEE 118-bus comparisons show that this modeling choice changes the
predicted pattern of secondary exposure: the initiating line need not
be the monitored element with sustained overload, and the stochastic
and deterministic descriptions identify different extents of exposure
within the tested ensemble.

The nominal probability law selects a line uniformly and a continuous
fault duration with mean $0.1\,\mathrm{s}$. CEM changes the proposal,
not that target, and strict exceedance is used throughout. The displayed
sampling runs have lower recorded CEM costs at their selected budgets,
but neither uniform variance reduction nor a real-time operational
update follows from those tables. Per-trajectory calculations and total
ensemble runtimes are distinct, and the probability-summary limitations
in Sec.~\ref{BenchAn} remain relevant to any precision claim.

The main modeling conclusion is that the prescribed variability matters
to the reduced predictions. Establishing its physical magnitude and
appropriate representation requires a reference-based study. In
particular, an extension should identify effective parameters from
high-fidelity electromechanical observations and examine how those
parameters depend on the pre-fault operating state. That extension may
modify the diffusion structure and stochastic convention; it should not
be described as fitting an unchanged model solely because both
representations are linear SDEs with fault-local noise.

Further work should retain the distinction between numerical integration
accuracy, sampling uncertainty, and physical-model discrepancy.
Voltage-coupled dynamics and protection-specific switching sequences are
needed where they materially change the screening observable. The
framework is intended to prioritize high-fidelity follow-up, rather than
replace EMT simulation, infer thermal damage from an exposure threshold,
or prescribe protection actions from an illustrative risk category.

\appendices
\numberwithin{equation}{section}
\section{Fault-Local Phase Noise: Flow, Sampling, and Moments}
\label{app:localized_noise_foundation}
This appendix develops the solution and moment equations for
\eqref{eq:swing-lin-matrix}, conditional on the faulted line,
duration $\tau$, operating condition, and initial state. The
interpretation remains It\^o. For compactness, boldface is suppressed
below: $x=\bm x$, $P=\bm P$, $A_0=\bm A_0$, and
$A_{\rm f}=\widetilde{\bm A}=A_0+\delta A$. Write $G=G_{ij}$,
$d=2n$, and $L_{\rm f}=\overline L(t)$ on $0\leq t<\tau$.
The positive diagonal matrix $M$ acts on the retained differential
coordinates. No commutation assumption is imposed on $A_{\rm f}$ and $G$.

\subsection{Meaning and structure of the diffusion}
\label{app:noise_structure}
Let $v=e_i-e_j$ and $B=\sigma_{ij}M^{-1}vv^{\mathsf T}$.
The lower-right diffusion block gives
\begin{align}
 d\omega&=(-M^{-1}D\omega-M^{-1}\overline L(t)\theta
                +M^{-1}p)\,dt,\nonumber\\
 d\theta&=\omega\,dt+
          \mathbb I\{t\leq\tau\}B\theta\,dW_t.
 \label{eq:noiseapp:block_sde}
\end{align}
Thus $\omega$ is not the pathwise derivative of the noisy phase.
A literal stochastic perturbation of the line Laplacian would enter
$ d\omega$, not $d\theta$. The specified diffusion is an effective
phase correction supported at the fault endpoints. Its two nonzero
components share one Wiener driver, and its coefficient $B$ has units
$\mathrm{s}^{-1/2}$. This interpretation does not identify the correction
with a measured susceptance process.

With $a=(0,M^{-1}v)^{\mathsf T}$ and $b=(0,v)^{\mathsf T}$, define
\begin{equation}
 c=b^{\mathsf T}a>0,\qquad
 \Pi=\frac{ab^{\mathsf T}}c,\qquad \chi_{ij}=\sigma_{ij}c.
 \label{eq:noiseapp:projector}
\end{equation}
Then $\Pi^2=\Pi$, $G=\chi_{ij}\Pi$, and $G^2=\chi_{ij} G$.
The matrix $G$ is not nilpotent for nonzero $\sigma_{ij}$.
For $K_{\rm f}=M^{-1}L_{\rm f}$,
\begin{equation}
 [A_{\rm f},G]=
 \begin{bmatrix}0&-K_{\rm f}B\\-B&0\end{bmatrix}\ne0
 \quad\text{if }\sigma_{ij}\ne0.
 \label{eq:noiseapp:commutator}
\end{equation}
Consequently, the usual single-exponential formula for commuting linear
SDEs does not give the coupled flow here
\cite{r4a}.

\subsection{Exact flow and reclosure}
\label{app:noise_flow}
The coefficients are globally Lipschitz with linear growth, giving a
unique continuous strong solution and finite moments on every finite
horizon for the specified initial state
\cite{noiseKaratzas1991}.
On $0\leq s\leq t\leq\tau$, define
\begin{equation}
 d\Phi(t,s)=A_{\rm f}\Phi(t,s)\,dt+G\Phi(t,s)\,dW_t,
 \qquad\Phi(s,s)=I_d.
 \label{eq:noiseapp:flow}
\end{equation}
The It\^o product rule applied to $\Phi(t,0)^{-1}x_t$ gives
\begin{equation}
 x_t=\Phi(t,0)x_0+\int_0^t\Phi(t,s)P\,ds.
 \label{eq:noiseapp:variation}
\end{equation}
The last integral is an ordinary time integral along the stochastic
flow \cite{r4b}.
This is an exact representation, not a claim that $\Phi$ depends only
on the final Brownian increment. Equivalently,
\begin{equation}
 x_t=x_{\rm det}(t)+
       \int_0^t e^{A_{\rm f}(t-s)}Gx_s\,dW_s.
 \label{eq:noiseapp:mild}
\end{equation}
Here $x_{\rm det}(t)$ is the deterministic response with the same
$A_{\rm f}$, $P$, and $x_0$, obtained by setting $G=0$ in
\eqref{eq:noiseapp:variation}. The random integrand prevents replacing
this convolution by a deterministic-kernel Gaussian draw.

The rank-one structure also gives the exact integrating factor
\begin{align}
 R_t&=I_d+
  \left(e^{\chi_{ij} W_t-\chi_{ij}^2t/2}-1\right)\Pi,\nonumber\\
 R_t^{-1}&=I_d+
  \left(e^{-\chi_{ij} W_t+\chi_{ij}^2t/2}-1\right)\Pi.
 \label{eq:noiseapp:factor}
\end{align}
Since $GP=0$, setting $x_t=R_ty_t$ reduces the SDE to
\begin{equation}
 \dot y_t=R_t^{-1}A_{\rm f}R_ty_t+P,\qquad y_0=x_0.
 \label{eq:noiseapp:random_ode}
\end{equation}
This is an ordinary differential equation for each continuous Brownian
path. Its coefficient remains random and time dependent; no commuting
matrix approximation has been made.

At reclosure, the state remains continuous. For $t>\tau$,
\begin{equation}
 x_t=e^{A_0(t-\tau)}x_\tau+
          \int_\tau^t e^{A_0(t-s)}P\,ds.
 \label{eq:noiseapp:post}
\end{equation}
The continuation is deterministic conditional on $x_\tau$, but its
unconditional law is still random. No reset to the mean is justified.

\subsection{Time-resolved sampling}
\label{app:noise_sampling}
The reported stochastic simulations use scalar-driver Milstein.
For a step wholly inside the active interval, its update is
\begin{align}
 x_{k+1}={}&x_k+(A_{\rm f}x_k+P)h_k+Gx_k\Delta W_k\nonumber\\
 &+\tfrac12G^2x_k\bigl((\Delta W_k)^2-h_k\bigr),
 \qquad \Delta W_k\sim\mathcal N(0,h_k).
 \label{eq:noiseapp:milstein}
\end{align}
The same scalar increment acts at both endpoints. A grid containing
$\tau$ avoids mixing the two coefficient regimes. For these smooth
linear coefficients, this construction has the standard strong-order-one
finite-horizon convergence under refinement
\cite{noiseKloeden1992,noiseHigham2001}.
After $\tau$, use~\eqref{eq:noiseapp:post}; no new Wiener increments
are drawn.

For comparison at the level of numerical constructions, an exponential
splitting alternative first applies the exact noise substep
\begin{equation}
 x^+=\left[I_d+
 \left(e^{\chi_{ij}\Delta W_k-\chi_{ij}^2h_k/2}-1\right)\Pi\right]x_k
 \label{eq:noiseapp:noise_step}
\end{equation}
and then the affine deterministic flow over $h_k$.
The individual substeps are exact, but their composition approximates
the coupled SDE because of~\eqref{eq:noiseapp:commutator}.
This alternative is not the production sampler used for the reported
stochastic experiments.
Neither construction is an exact Gaussian terminal-state sampler.
A path-dependent overload calculation requires temporally consistent
samples and a separate refinement of its observation grid.

\subsection{Exact first and second moments}
\label{app:noise_moments}
Let $\mu=\mathbb E[x]$, $R_2=\mathbb E[xx^{\mathsf T}]$, and
$C=R_2-\mu\mu^{\mathsf T}$. The It\^o product rule yields
\begin{align}
 \dot\mu&=\overline A(t)\mu+P,\label{eq:noiseapp:mean}\\
 \dot R_2&=\overline A(t)R_2+R_2\overline A(t)^{\mathsf T}
              +P\mu^{\mathsf T}+\mu P^{\mathsf T}\nonumber\\
 &\quad+\mathbb I\{t\leq\tau\}GR_2G^{\mathsf T},
 \label{eq:noiseapp:second}\\
 \dot C&=\overline A(t)C+C\overline A(t)^{\mathsf T}\nonumber\\
 &\quad+\mathbb I\{t\leq\tau\}
          G(C+\mu\mu^{\mathsf T})G^{\mathsf T}.
 \label{eq:noiseapp:cov}
\end{align}
Initially $\mu(0)=x_0$ and $C(0)=0$. These equations close exactly;
no Gaussian assumption is needed
\cite{noiseJimenez2015,noiseJimenezOzaki2002}.
The fixed-drift It\^o mean does not depend on $\sigma_{ij}$.
For $t>\tau$,
\begin{equation}
 C(t)=e^{A_0(t-\tau)}C(\tau)e^{A_0^{\mathsf T}(t-\tau)}.
 \label{eq:noiseapp:cov_post}
\end{equation}

For one constant segment with drift $A$ and diffusion $H$, put
\begin{align}
 \mathcal K_{A,H}&=I_d\otimes A+A\otimes I_d+H\otimes H,\nonumber\\
 J_P&=I_d\otimes P+P\otimes I_d,\nonumber\\
 z&=(\operatorname{vec}R_2,\mu,1)^{\mathsf T},\nonumber\\
 \mathcal T_{A,H}&=
 \begin{bmatrix}\mathcal K_{A,H}&J_P&0\\0&A&P\\0&0&0\end{bmatrix}.
 \label{eq:noiseapp:lift}
\end{align}
With column-wise vectorization,
\begin{equation}
 z(t+h)=e^{h\mathcal T_{A,H}}z(t).
 \label{eq:noiseapp:exact_moments}
\end{equation}
Use $(A,H)=(A_{\rm f},G)$ before reclosure and $(A_0,0)$ afterward.
This gives exact continuous moment propagation up to numerical
exponential evaluation. A matrix-free action can avoid forming the
full Kronecker matrix, although storing a full covariance remains
quadratic in $d$. Freezing $G\mu(s)\mu(s)^{\mathsf T}G^{\mathsf T}$
at a step endpoint is instead a time-discretization approximation.

For an ensemble of independently sampled scenarios $Z$, conditional
moments must be combined as
\begin{align}
 \mu_{\rm mix}&=\mathbb E_Z[\mu_Z],\nonumber\\
 C_{\rm mix}&=\mathbb E_Z[C_Z]+\operatorname{Cov}_Z(\mu_Z).
 \label{eq:noiseapp:scenario_moments}
\end{align}
This distinguishes within-scenario Wiener variability from variability
across fault locations, durations, or operating states.

The equivalent Stratonovich drift of the present It\^o equation is
$\overline A(t)-\tfrac12\mathbb I\{t\leq\tau\}G^2$.
A model instead postulated with Stratonovich drift $\overline A(t)$ has
a different It\^o drift and can have a noise-dependent mean.
Moment formulas must use the drift corresponding to the stated
convention. Mathematical well-posedness, numerical convergence, and
physical validation remain distinct questions.

\section{Milstein Implementation and Algebraic Verification}
\label{app:p2_milstein_analysis}
The stochastic simulations use Milstein. This appendix makes its
matrix and component forms explicit and derives its exact discrete
moments. The derivation complements the continuous moment equations
in Appendix~\ref{app:noise_moments}; it is not a claim of a newly
executed convergence experiment.

\subsection{Endpoint update and switching}
Use the notation of Appendix~\ref{app:localized_noise_foundation}.
Conditional on a sampled $(e,\tau)$ and a specified time grid,
$\Delta W_k=\sqrt{h_k}\,\xi_k$, where the $\xi_k$ are independent
standard normal variables. Each trajectory has its own sequence.
The same scalar $\Delta W_k$ acts at both endpoints of its faulted
line. Independent endpoint increments would define a different SDE.
Since $G^2=\chi_{ij}G$, the fault-on step becomes
\begin{align}
 \omega_{k+1}&=\omega_k+h_k M^{-1}
                 (p-D\omega_k-L_{\rm f}\theta_k),\nonumber\\
 \theta_{k+1}&=\theta_k+h_k\omega_k+B\theta_k\eta_k,\nonumber\\
 \eta_k&=\Delta W_k+\frac{\chi_{ij}}2
                         ((\Delta W_k)^2-h_k),\nonumber\\
 B&=\sigma_{ij}M^{-1}v_{ij}v_{ij}^{\mathsf T}.
 \label{eq:p2_milstein_components}
\end{align}
All right-hand sides use the state at $t_k$. The diffusion action is
\begin{equation}
 (B\theta)_r=\begin{cases}
 \sigma_{ij}(\theta_i-\theta_j)/m_i,&r=i,\\
 -\sigma_{ij}(\theta_i-\theta_j)/m_j,&r=j,\\
 0,&r\notin\{i,j\}.
 \end{cases}
 \label{eq:p2_endpoint_action}
\end{equation}
Its evaluation requires the endpoint difference, not a dense product.
The symbol $\chi_{ij}$ is the diffusion coefficient from
\eqref{eq:noiseapp:projector}; it is distinct from the dimensionless
CEM retention factor $\kappa=0.9$.

The integration grid contains $\tau$ whenever $0<\tau<T$.
A step ending at reclosure uses the fault-on coefficients; the following
segment uses $A_0$ and zero diffusion. The state reached at $\tau$ is
retained for the deterministic continuation
\eqref{eq:noiseapp:post}. For $\tau\geq T$, only the fault-on branch
is evaluated through $T$. None of these rules clips the sampled
$\tau$ in its probability weight.

Milstein has strong order one for these globally Lipschitz linear
coefficients on a fixed finite horizon and an event-aligned grid
\cite{noiseKloeden1992,noiseHigham2001}. This is a convergence statement
under refinement, not a numerical error bound at an unspecified step
size. The run-specific step sizes, random seeds, and observation grid
are not present in the retained result tables.

\subsection{Exact moments of the discretized method}
Let $F_k=I+h_kA_{\rm f}$, $b_k=h_kP$, and
$R_k=\mathbb E[x_kx_k^{\mathsf T}]$. Write the random step matrix as
\begin{equation}
\begin{aligned}
 M_k&=F_k+G\Delta W_k+\tfrac12G^2((\Delta W_k)^2-h_k),\\
 x_{k+1}&=M_kx_k+b_k.
\end{aligned}
 \label{eq:p2_milstein_random_matrix}
\end{equation}
For a grid fixed conditional on the initiating scenario, the new
increment is independent of $x_k$. Its centered moments satisfy
\begin{align}
 \mathbb E[\Delta W_k]&=0,\qquad
 \mathbb E[(\Delta W_k)^2]=h_k,\nonumber\\
 \mathbb E[\Delta W_k((\Delta W_k)^2-h_k)]&=0,\nonumber\\
 \mathbb E[((\Delta W_k)^2-h_k)^2]&=2h_k^2.
 \label{eq:p2_increment_identities}
\end{align}
Expanding \eqref{eq:p2_milstein_random_matrix} therefore gives
\begin{align}
 \mu_{k+1}&=F_k\mu_k+b_k,\label{eq:p2_milstein_mean}\\
 R_{k+1}&=F_kR_kF_k^{\mathsf T}
       +(F_k\mu_k)b_k^{\mathsf T}
       +b_k(F_k\mu_k)^{\mathsf T}+b_kb_k^{\mathsf T}\nonumber\\
 &\quad+h_kGR_kG^{\mathsf T}
       +\frac{h_k^2}{2}G^2R_k(G^2)^{\mathsf T}.
 \label{eq:p2_milstein_second}
\end{align}
Subtracting $\mu_{k+1}\mu_{k+1}^{\mathsf T}$ yields
\begin{align}
 C_{k+1}&=F_kC_kF_k^{\mathsf T}
       +h_kG(C_k+\mu_k\mu_k^{\mathsf T})G^{\mathsf T}\nonumber\\
 &\quad+\frac{h_k^2}{2}G^2(C_k+\mu_k\mu_k^{\mathsf T})(G^2)^{\mathsf T}.
 \label{eq:p2_milstein_covariance}
\end{align}
These are exact expectations of the discrete Milstein chain, not an
approximation obtained by freezing the continuous covariance forcing.
They require no Monte Carlo sampling and retain the nonzero $G^2$
correction. Each covariance contribution is positive semidefinite when
$C_k$ is positive semidefinite, in exact arithmetic.

The discrete mean follows explicit Euler during the active interval;
it is not the exact continuous deterministic trajectory at finite
$h_k$. For a fixed incoming mean and constant coefficients, the
one-step difference is
\begin{equation}
 \mu_{\rm exact}(t_k+h_k)-(F_k\mu_k+b_k)
 =\frac{h_k^2}{2}(A_{\rm f}^2\mu_k+A_{\rm f}P)
     +\mathcal O(h_k^3).
 \label{eq:p2_mean_step_defect}
\end{equation}
The zero-noise limit is consequently the Euler discretization during
that interval, followed by the stated deterministic continuation, not
an exact finite-step solution of the fault-on ODE. Under refinement,
the discrete moments approach the continuous equations in
Appendix~\ref{app:noise_moments}. A numerical comparison of those two
systems isolates discretization bias; comparison of sampled moments
with \eqref{eq:p2_milstein_mean}--\eqref{eq:p2_milstein_covariance}
checks the random implementation at a fixed grid.

For a coupled refinement, a coarse increment is the sum of the fine
increments covering the same interval. Drawing an independent coarse
path measures path variability as well as numerical error. The
analytical identities above supply verification references, but no
sampled-moment or coupled-refinement results for the published IEEE
118-bus runs are available in the retained record. Such tests are not
reported as completed here.

\subsection{Noise units and numerical normalization}
Angles are dimensionless and time is measured in seconds in the
continuous equations. If power is measured in a common unit,
$[M]=[p]\,\mathrm{s}^2$, $[D]=[p]\,\mathrm{s}$, and
$[\beta]=[p]$. The phase diffusion requires
\begin{equation}
 [B]=\mathrm{s}^{-1/2},\qquad
 [\sigma_{ij}]=[M]\,\mathrm{s}^{-1/2}
                =[p]\,\mathrm{s}^{3/2}.
 \label{eq:p2_noise_units}
\end{equation}
Thus the reported numerical assignment $\sigma_{ij}=\beta_{ij}$ is
not a dimensional equality.

An explicit change of numerical bases makes the distinction precise.
For power scale $S_\star$ and time scale $t_\star$, define
\begin{align}
 t'&=t/t_\star,\qquad W'_{t'}=W_{t_\star t'}/\sqrt{t_\star},\nonumber\\
 \omega'&=t_\star\omega,\qquad
 M'=M/(S_\star t_\star^2),\qquad D'=D/(S_\star t_\star),\nonumber\\
 L'&=L/S_\star,\qquad p'=p/S_\star,\qquad
 \sigma'_{ij}=\sigma_{ij}/(S_\star t_\star^{3/2}).
 \label{eq:p2_nondimensionalization}
\end{align}
The primed SDE has the same block structure and standard Wiener driver.
Choosing stored coefficients $\sigma'_{ij}=\beta'_{ij}$ prescribes a
model in those bases; its dimensional counterpart is
$\sigma_{ij}=t_\star^{3/2}\beta_{ij}$. The numerical power and inertia
bases used in the historical realization are not specified in the
retained experiment record, so no physical noise amplitude is inferred
from the stored equality. The lower-right block itself remains
unchanged.

Unlike independent absolute-endpoint corrections, the line-difference
matrix here annihilates a uniform phase shift:
$B\boldsymbol1=0$. Changing the initial phase reference by a constant
therefore leaves both the drift's line differences and this diffusion
unchanged. It does not make $\omega$ the ordinary pathwise derivative
of the noisy phase.

\section{Sampling and Experimental Specification}
\label{app:p2_experiment_record}
\subsection{Catalogue, fixed budget, and timing}
The supplied CEM configuration is
$N_{\rm CE}=200$, $K=15$, $n_{\rm elite}=20$, $\rho=0.1$,
$\kappa=0.9$, $E=177$, $\phi_e^{(0)}=1/177$, and
$r^{(0)}=1.0\,\mathrm{s}$. Equations
\eqref{eq:categorical_exponential_proposal}--\eqref{eq:proposal_support_update}
retain the line-indexed family and its shared-mean restriction. The
scalar $r^{(0)}$ initializes the duration mean but does not identify
its tying rule after adaptation. The final proposal record is not
available. The nominal mean is fixed at
$0.1\,\mathrm{s}$ in both CEM and direct MC.

The catalogue contains the $177$ candidate transmission lines, so
$\pi_e=\phi_e^{(0)}=1/177$ and
$\sum_{e\in\mathcal F}\pi_e=177(1/177)=1$.
The nominal law and likelihood ratio use the same $E=177$ as the
initialization. After adaptation, the likelihood-ratio denominator
retains the learned $\phi_e^\ast$ rather than replacing it by $1/177$.
The separately reported transformers can be monitored elements without
being additional entries in the initiating-fault catalogue.
The model's differential matrices require the actual inertia and
damping assignments on the retained coordinates. The power-flow case
reference alone does not specify those dynamic coefficients.

The adaptation and final-batch counts satisfy
\begin{equation}
 N_{\rm adapt}=200\times15=3{,}000,\qquad
 N_{\rm total}=N+3{,}000.
 \label{eq:p2_evaluation_budget}
\end{equation}
For the final sizes in Tables~\ref{table4}--\ref{table5}, this gives
\begin{center}
\begin{tabular}{r|rrr}
\toprule
$\gamma\ [\mathrm{s}]$ & $N$ & $N_{\rm adapt}$ & $N_{\rm total}$\\
\midrule
$0$ & $2{,}000$ & $3{,}000$ & $5{,}000$\\
$0.5$ & $3{,}500$ & $3{,}000$ & $6{,}500$\\
$5$ & $10{,}000$ & $3{,}000$ & $13{,}000$\\
$10$ & $12{,}500$ & $3{,}000$ & $15{,}500$\\
\bottomrule
\end{tabular}
\end{center}
These are derived evaluation counts for the supplied schedule and one
final batch per threshold. They are not new timing results. The
original timing values are preserved, including the differences between
the stochastic and deterministic tables. Their repetition counts,
hardware, software versions, and inclusion of preprocessing or
adaptation are not specified. A stage-by-stage runtime cannot be
recovered from those means alone.

The positive MC summaries $9.85\times10^{-6}$ and
$1.89\times10^{-6}$ at $N=50{,}000$ cannot be single-batch event
frequencies: the smallest positive such frequency is
$2\times10^{-5}$. Neither value can be explained by rounding that
minimum to the displayed precision. The available tables do not
identify a repeated-batch or alternative estimator construction.
Those entries remain reported summaries and are excluded from claims
of verified MC tail resolution or estimator convergence.

\subsection{Conditioning and duration-resolved plots}
\label{app:p2_duration_score}
For an initiating category $e$ and fixed duration $u$, define
\begin{equation}
\begin{aligned}
 q_b(e,u)&=\mathbb P_W[S_b(x^{e,u,W})>T^{(*)}],\\
 q_b(u)&=\frac1E\sum_{e\in\mathcal F}q_b(e,u).
\end{aligned}
 \label{eq:p2_fixed_duration_risk}
\end{equation}
The first quantity averages only over the Wiener path; the second also
averages over initiating categories. In the deterministic comparator,
$q_b(e,u)$ is an indicator for fixed initial conditions. The nominal
probability integrates over duration:
\begin{equation}
 Q_b(T^{(*)})=\int_0^\infty q_b(u)
                    \lambda_0^{-1}e^{-u/\lambda_0}\,du.
 \label{eq:p2_duration_marginal}
\end{equation}
Neither $\mathbb P[S_b>T^{(*)},\tau\leq u]$ nor
$\mathbb P[S_b>T^{(*)}\mid\tau\leq u]$ equals $q_b(u)$ in general.
Their numerical denominators are also different. The graphical
summaries in Figs.~\ref{fig3} and~\ref{fig5} do not identify which
conditional or cumulative estimator was used, its duration bins, or
the sample counts and weights at each displayed duration.

A nominal batch of $100{,}000$ exponential durations with mean
$0.1\,\mathrm{s}$ has expected counts of approximately $33.55$
beyond $0.8\,\mathrm{s}$, $4.540$ beyond $1\,\mathrm{s}$, and
$0.03059$ beyond $1.5\,\mathrm{s}$. These follow from
$N\exp(-u/\lambda_0)$ and are expectations, not recovered counts.
The reported total budget therefore does not document dense
long-duration conditional curves. Their original generation and
plotting records are required before assigning a specific estimator to
them. The exposure thresholds and the confirmed nominal law are not
changed to accommodate the graphics.

\subsection{Trajectory grid, observation grid, and strict scoring}
Milstein advances the state on an integration grid. The overload
observable is evaluated on an observation grid, which need not be the
same grid. For illustration of the distinction, the explicitly defined
left-endpoint score on $0=s_0<\cdots<s_J=T$ is
\begin{equation}
 S_{b,h,\Delta}^{\rm left}
 =\sum_{j=0}^{J-1}(s_{j+1}-s_j)
     \mathbb I\{|\beta_b(s_j)v_b^{\mathsf T}\theta_h(s_j)|>\bar p_b\}.
 \label{eq:p2_discrete_score}
\end{equation}
The coupling here is the deterministic effective coupling on the
appropriate side of reclosure, not a realization of pointwise white
noise. Equality to the flow threshold does not contribute, and the
final event also uses $S_{b,h,\Delta}>T^{(*)}$ or
$S_{h,\Delta}>\gamma$ strictly. The formula defines a quadrature
choice; the retained experiments do not identify their quadrature
rule, internal step size, observation interval, or interpolation.
No numerical values are assigned to those missing settings.

Importance weighting preserves the probability of the implemented
score, conditional on the fitted proposal, with the physical and numerical
scoring specification fixed:
\begin{equation}
 \mathbb E_q[w\,\mathbb I\{S_{h,\Delta}>\gamma\}
                       \mid\mathcal F_{\rm adapt}]
     =\mathbb P_{p_Z,W}[S_{h,\Delta}>\gamma].
 \label{eq:p2_discrete_target}
\end{equation}
Here $\mathcal F_{\rm adapt}$ denotes the adaptation history.
This identity does not establish accuracy relative to the continuous
score. For trajectories evaluated on a common probability space,
\begin{equation}
 |Q_{h,\Delta}(\gamma)-Q(\gamma)|
 \leq\mathbb P\!\left[|S-\gamma|\leq|S-S_{h,\Delta}|\right].
 \label{eq:p2_threshold_error}
\end{equation}
A difference between the strict event indicators requires the score
error to span the threshold. Consequently, a state-convergence order
does not automatically give the same convergence rate for the
discontinuous overload event. Observation refinement and trajectory
refinement address different parts of that error.

\subsection{Scope of the stochastic--deterministic comparison}
The confirmed common specifications are the case designation,
initial operating-state description, two-thirds coupling during the
modified interval, exponential duration mean, and strict
$T^{(*)}=1\,\mathrm{s}$ exposure threshold. The stochastic calculation
adds the stated phase diffusion and uses Milstein. No record establishes
paired fault draws, identical integration errors, or a complete
selection rule for every displayed monitored line. The reported
$45$-versus-$8$ count is therefore a comparison of the retained
selections, not a false-positive rate, false-negative rate, or full
population classification statistic.

An exact reconstruction also requires the indexed case realization,
retained state dimension, numerical $M$, $D$, $p$, coupling and flow-limit
arrays, numerical bases, source revision, seeds, final frozen proposal,
and run-to-figure association. The mathematical specification and the
additional supplied settings do not provide those empirical records.
The analytical moment identities in this paper are established from
the stated SDE and discretization; historical numerical validation,
physical calibration, and complete experimental reproducibility remain
separate claims.

\bibliographystyle{ieeetr}
\bibliography{refs,Noise_References}
\end{document}